\documentclass[12pt]{article}
\usepackage{amsfonts,amssymb,amsmath,eepic,tikz}
\usetikzlibrary{trees}

\newenvironment{tikz-trees-diagram}
{                                                                                                                                                        \begin{tikzpicture}
  [edge from parent fork down,
   level distance = 0.5in,
   sibling distance = 0.5in,
   inner sep = 1mm,
   scone/.style = {draw=blue!50, fill=blue!5, rectangle,
     thick},
   nscone/.style = {draw=red!70, fill=red!5, double, rectangle,
     thick, inner sep=2mm},
   leaf/.style = {draw=black, fill=white, circle, thick},
   every label/.style = {red}]
}
{\end{tikzpicture}}

\usepackage{booktabs}

 \def\medskipamount{12pt} \def\smallskipamount{6pt}

\newcounter{bitcount}
\newcommand{\bit}[1]{\addtocounter{bitcount}{1}\pagebreak[3]
\subsection{#1}\nopagebreak\setcounter{equation}{0}}

\renewcommand{\theequation}{\thesubsection .\arabic{equation}}
\renewcommand{\thesubsection}{\arabic{bitcount}}

\newcommand{\re}[1]{\mbox{\bf (\ref{#1})}}

\catcode`\@=\active
\catcode`\@=11
\def\@eqnnum{\hbox to .01pt{}\rlap{\bf \hskip -\displaywidth(\theequation)}}
\catcode`\@=12

\begin{document}


\catcode`\@=\active
\catcode`\@=11
\newcommand{\nc}{\newcommand}


\nc{\bs}[1]{ \addvspace{\medskipamount} \pagebreak[3]
\refstepcounter{equation}
\noindent {\bf (\theequation) #1.} \em \nopagebreak}

\nc{\es}{\rm \par \addvspace{\medskipamount} } 

\nc{\ess}{\end{em} \par}

\nc{\br}[1]{ \addvspace{\medskipamount} \pagebreak[3]
\refstepcounter{equation} 
\noindent {\bf (\theequation) #1.} \nopagebreak}

\nc{\brs}[1]{ \pagebreak[3]
\refstepcounter{equation} 
\noindent {\bf (\theequation) #1.} \nopagebreak}

\nc{\er}{\par \addvspace{\medskipamount} }


\nc{\vars}[2]
{{\mathchoice{\mb{#1}}{\mb{#1}}{\mb{#2}}{\mb{#2}}}}
\nc{\A}{\mathbb A}
\nc{\C}{\mathbb C}
\renewcommand{\H}{H}
\nc{\N}{\mathbb N}
\renewcommand{\O}{{\mathcal O}}
\nc{\Q}{\mathbb Q}
\nc{\R}{\mathbb R}
\nc{\Z}{\mathbb Z}
\renewcommand{\P}{\mathbb P}


\nc{\oper}[1]{\mathop{\mathchoice{\mbox{\rm #1}}{\mbox{\rm #1}}
{\mbox{\rm \scriptsize #1}}{\mbox{\rm \tiny #1}}}\nolimits}
\nc{\Adj}{\oper{Adj}}
\nc{\Hull}{\oper{Hull}}
\nc{\diag}{\oper{diag}}
\nc{\id}{\oper{id}}
\nc{\Proj}{\oper{Proj}}
\nc{\Spec}{\oper{Spec}}
\nc{\Gr}{\oper{Gr}}

\nc{\operlim}[1]{\mathop{\mathchoice{\mbox{\rm #1}}{\mbox{\rm #1}}
{\mbox{\rm \scriptsize #1}}{\mbox{\rm \tiny #1}}}}

\nc{\lp}{\raisebox{-.1ex}{\rm\large(}}
\nc{\rp}{\raisebox{-.1ex}{\rm\large)}}

\nc{\al}{\alpha}
\nc{\be}{\beta}
\nc{\la}{\lambda}
\nc{\La}{\Lambda}
\nc{\ep}{\varepsilon}
\nc{\si}{\sigma}
\nc{\om}{\omega}
\nc{\Om}{\Omega}
\nc{\Ga}{\Gamma}
\nc{\Si}{\Sigma}


\nc{\Left}[1]{\hbox{$\left#1\vbox to
    11.5pt{}\right.\nulldelimiterspace=0pt \mathsurround=0pt$}}
\nc{\Right}[1]{\hbox{$\left.\vbox to
    11.5pt{}\right#1\nulldelimiterspace=0pt \mathsurround=0pt$}}

\nc{\updown}{\hbox{$\left\updownarrow\vbox to
    10pt{}\right.\nulldelimiterspace=0pt \mathsurround=0pt$}}


\nc{\beqas}{\begin{eqnarray*}}
\nc{\co}{{\cal O}}
\nc{\cx}{{\C^{\times}}}
\nc{\down}{\Big\downarrow}
\nc{\Down}{\left\downarrow
    \rule{0em}{8.5ex}\right.}
\nc{\downarg}[1]{{\phantom{\scriptstyle #1}\Big\downarrow
    \raisebox{.4ex}{$\scriptstyle #1$}}}
\nc{\eeqas}{\end{eqnarray*}}
\nc{\fp}{\mbox{     $\square$} \par \addvspace{\smallskipamount}}
\nc{\lrow}{\longrightarrow}
\nc{\pf}{\noindent {\em Proof}}
\nc{\sans}{\, \backslash \,}
\nc{\st}{\, | \,}

\nc{\beq}{\begin{equation}}
\nc{\eeq}{\end{equation}}

\nc{\GL}[1]{\gp{GL}{#1}}

\hyphenation{para-met-riz-ing sub-bundle}

\catcode`\@=12



\noindent
{\LARGE \bf Resolving toric varieties with Nash blow-ups}
\medskip \\ 
{\bf Atanas Atanasov, Christopher Lopez, Alexander Perry} \\ 
Department of Mathematics, Columbia University \\
2990 Broadway, New York, N.Y. 10027 \smallskip\\
{\bf Nicholas Proudfoot}\\
Department of Mathematics, University of Oregon\\
1222 University of Oregon, Eugene, Ore. 97403 \smallskip \\
{\bf Michael Thaddeus } \\ 
Department of Mathematics, Columbia University \\
2990 Broadway, New York, N.Y. 10027
\renewcommand{\thefootnote}{}
\footnotetext{N.P. supported by NSF grant
  0738335;  M.T. supported by NSF grants 0401128 and 0700419.} 

\bit{Introduction} 

Let $X$ be a variety over an algebraically closed field $K$.
Its {\it Nash blow-up} is a variety over $K$ with a projective
morphism to $X$, which is an isomorphism over the smooth locus.
Roughly speaking, it parametrizes all limits of tangent planes
to $X$ (a precise definition is given in \S2 below).  The Nash
blow-up of a singular $X$ is not always smooth but seems, in
some sense, to be less singular than $X$.  Strictly speaking this is
false, for in characteristic $p>0$, as explained by Nobile
\cite{nobile}, the plane curve $x^p - y^q = 0$ is its own Nash
blow-up for any $q>0$.  In this and other ways the ordinary
Nash blow-up proves intractable.

However, let the {\it normalized Nash blow-up} be the
normalization of the Nash blow-up.  Then, of course, the
normalized Nash blow-up of every curve is smooth.  The
normalized Nash blow-up of a surface can be singular, but
Hironaka \cite{hironaka} and Spivakovsky \cite{bourbaki,
  spivakovsky} have shown that every surface becomes smooth
after finitely many normalized Nash blow-ups.  Thus we are drawn to ask
the following.

\bs{Question} Is every variety desingularized by 
finitely many normalized Nash blow-ups?
\es

According to Spivakovsky \cite{spivakovsky}, Nash asked
Hironaka this question in the early 1960s.  An affirmative
answer would give a canonical procedure for desingularizing an
arbitrary variety.  The answer to this question is not known
and is surely difficult.  In this paper we address a more narrow
question.

\bs{Question} Is every toric variety desingularized by finitely
many normalized Nash blow-ups?
\label{query}
\es

We do not answer this question conclusively either.  But we do
exhibit abundant evidence supporting an affirmative answer.
Using the ``toric dictionary,'' which translates every problem in
toric geometry into a problem on convex polyhedra, we 
convert Question \re{query} into a problem amenable to computer
calculation.  Then we carry out this calculation for over a thousand
examples.  In every case, finitely many Nash blow-ups
produce a smooth toric variety.

{\it Summary of the paper.}  In \S2 we introduce the
toric dictionary and, following Gonzalez-Sprinberg \cite{thesis},
translate the question into convex geometry.  In \S3 we
summarize the effect of the iterated Nash blow-up, in the toric
case, using the notion of a {\it resolution
  tree}.  In \S4 we spell out what happens in the
2-dimensional toric case in terms of continued fractions.  In
\S5 we digress briefly on the classification of quasi-smooth
affine toric varieties, those corresponding to simplicial cones
in the toric dictionary.  Then in \S6 we give an account of
our computer investigations.

{\it Notation and conventions.} We slightly abuse terminology in
two ways.  First, as we are concerned with the normalized
Nash blow-up throughout, we refer to it simply as the {\it Nash
blow-up}.  Second, as we are concerned with rational
polyhedra and rational polyhedral cones throughout, we refer to
them simply as {\it polyhedra\/} and {\it cones}.  We denote the
natural numbers, including 0, by $\Z_+$, and we likewise
denote the nonnegative rational numbers, including 0, by
$\Q_+$.  We denote the span of $v_1, \dots, v_k$ with
coefficients in $S$ by $S\langle v_1, \dots, v_k \rangle$.
Thus, for example, the first quadrant in $\Q^2$ is denoted
$\Q_+\langle e_1, e_2 \rangle$.

{\it Acknowledgements.}  The very helpful advice provided by Kevin
Purbhoo is gratefully acknowledged.  We also thank
Jeffrey Lagarias for an inspiring conversation, Dave Bayer for
recommending the use of {\tt 4ti2}, and Sam
Payne and Greg Smith for informing us of their parallel work on
the subject.

\bit{Equivalence to a combinatorial problem}

{\bf Nash blow-ups.}  Let $X \subset \P^n$ be a quasiprojective
variety of dimension $d$ over an algebraically closed field
$K$.  The {\it Gauss map} is the rational map $X \dasharrow
\Gr(d\!+\!1, \,n\!+\!1)$ taking a smooth point to its tangent
plane.  The {\it Nash blow-up} of $X$ is defined to be the
closure of the graph of the Gauss map.  The {\it normalized
  Nash blow-up} of $X$ is the normalization of the Nash blow-up
of $X$.  Gonzalez-Sprinberg's and Spivakovsky's results are
concerned with this variant, as are ours.  Consequently, we
shall abuse terminology by referring to a normalized Nash
blow-up simply as a Nash blow-up.

{\it Remarks:}

(1) As defined, the Nash blow-up appears to depend on the projective
embedding of $X$, but it can be reformulated in terms of
K\"ahler differentials and hence depends only on $X$ (and makes
sense even if $X$ is not quasiprojective) \cite{bourbaki}.

(2) Since the normalization of a variety over $K$ is a finite
morphism \cite[I 3.9A]{hartshorne} and the pullback of an ample
bundle by a finite morphism is ample \cite[1.7.7]{lazarsfeld},
the normalization is a projective morphism.  Hence the natural
morphism from the (normalized) Nash blow-up of $X$ to $X$ is
projective.

(3) Clearly the Nash blow-up of a smooth variety is itself, and
the Nash blow-up of a product is a product.

(4) If $X \subset \A^d$ is an affine variety, we may consider
the analogous construction using the Gauss map $X \dasharrow
\Gr(d,n)$, but this produces exactly the same thing, since the
morphism $X \times \Gr(d,n) \to X \times \Gr(d\!+\!1,
\,n\!+\!1)$ given by $(x,V) \mapsto \big( x, K\langle(1\times x) \rangle
\oplus (0\times V)\big)$ is a closed embedding.

\smallskip
\noindent
{\bf Toric varieties.}  We review here some standard
definitions and facts about toric varieties.  For proofs, we
refer the reader to Ewald \cite{ewald}, Fulton \cite{fulton},
Miller \& Sturmfels \cite{ms}, and Thaddeus \cite{thaddeus}.

A {\it polyhedron} in $\Q^d$ is a subset $P$ defined by
finitely many weak affine inequalities, say $\sum_{j=1}^d
a_{ij} x_j \geq b_i$.  It is a {\it polyhedral cone} if the
inequalities are linear, that is, all $b_i = 0$.  For
simplicity we refer to polyhedral cones simply as {\it cones}.
We also assume that all polyhedra are {\it rational}, meaning
that all $a_{ij}$ and all $b_i$ are rational.  A polyhedron is
{\it proper} if it contains no affine linear subspace besides a
point, and is contained in no affine linear subspace besides
$\Q^d$.

A {\it face} $F$ of $P$ is the locus where equality holds in
some fixed subset of the inequalities above.  It is a {\it
  facet\/} if its affine linear span has codimension 1 in that of
$P$.  It is a {\it vertex\/} if it is a point.  A proper cone in
$\Q^d$ is {\it simplicial\/} if it has exactly $d$ facets.  For any face $F
\subset P$, the {\it localization} $P_F$ is the cone generated,
as a semigroup, by the Minkowski difference $P - F$.

For $t \neq 0$, let $tP = \{ tv \st v \in P \}$; however, for
$t = 0$, by convention let $0P$ denote the {\it cone at
  infinity} defined by the same inequalities as $P$, except
with constant terms set to zero.  The reason for this
convention is that $\{(t,v) \in \Q_+ \times \Q^d \st v \in
tP\}$ is then a cone in $\Q^{d+1}$, defined by the inequalities
$\sum_{j=1}^d a_{ij} x_j \geq b_i x_0$ and $x_0 \geq 0$.  It is
called the {\it cone over $P$} and denoted $C(P)$.

A {\it torus} is a product of finitely many copies of the
multiplicative group of $K$.
A {\it toric variety} is a normal variety on which a torus acts with
finitely many orbits.  There is a one-to-one correspondence between
polyhedra with integer vertices and toric varieties that are
projective over an affine, equipped with a lifting of the torus action
to $\O(1)$.  It is given as follows.  For a polyhedron $P \subset
\Q^d$, the semigroup algebra $K[C(P) \cap \Z^{d+1}]$ is graded by the
0th coordinate.  Let $X(P) = \Proj K[C(P) \cap \Z^{d+1}]$.  This is a
quasiprojective variety acted on by the torus $T = \Spec K[\Z^d]$.
For example, if $P$ is already a cone, then $C(P) = \Q_+ \times P$ and
$X(P) = \Proj K[P \cap \Z^d] [x_0] = \Spec K[P \cap \Z^d]$, the affine
toric variety usually associated to a cone.  In general, $X(P)$ is
projective over the affine $X(0P)$, because $C(P) \cap (0 \times \Q^d)
= 0P$.

{\it Further remarks:}

(5) A polyhedron $P$ is proper if and only if (a) the torus action
on $X(P)$ is effective, and (b) $X(P)$ is not a direct product of a
toric variety with a torus.  So in light of remark (3), there
is no loss of generality, for the purposes of Nash blowing-up, in
assuming that $P$ is proper.

(6) Any toric variety has a natural cover by open affine toric
subvarieties.  Indeed, $X(P)$ is covered by the affine varieties
$X(P_F)$, where $F$ runs over the faces of $P$.  If $P$ is
proper, just the vertices are sufficient.

(7) Define two cones to be {\it equivalent\/} if an element of
$GL(d,\Z)$ takes one to the other.  Then equivalent cones
clearly lead to isomorphic toric varieties, with the torus action
adjusted by the appropriate automorphism of $T$.

(8) An affine toric variety $X(C)$, with $C$ proper, is smooth
if and only if it is isomorphic to $\A^d$, or equivalently, if
$C$ is equivalent to the orthant $\Q_+\langle e_1, \dots, e_d
\rangle$.

\smallskip
\noindent 
{\bf Nash blow-ups of toric varieties.}  Now let $C$
be a cone in $\Q^d$.  Let $H$ be the Hilbert basis of the
semigroup $C \cap \Z^d$, that is, the set of indecomposable
nonzero elements in the semigroup.  This is the unique minimal
set of generators of $C \cap \Z^d$.  By Gordan's lemma
\cite[\S1.2, Prop.\ 1]{fulton} $H$ is a finite set, say with $n$
elements.  Let $M$ be the $n \times d$ integer matrix whose
rows are the elements of $H$.

Let $S = \{h_1 + \cdots + h_d \, | \, h_i \in H \mbox{ linearly
  independent} \}$.  Since $S$ is finite, its convex hull is a
compact polyhedron $\Hull S$.  Hence the Minkowski sum $C +
\Hull S$ is a polyhedron whose cone at infinity is $C$.

The following result is proved (in the language of fans) by
Gonzalez-Sprinberg \cite{thesis}.

\bs{Theorem} 
\label{theorem:gonzalez}
The Nash blow-up of $X(C)$ is $X(C + \Hull S)$.
\es

\pf.  Without loss of generality $C$ may be assumed proper.  In
this case $X(C)$ has a unique $T$-fixed point $q$.

Let $X = X(C)$.  The Nash blow-up of $X$ is plainly a toric variety,
projective over $X$.  It is therefore $X(P)$ for some polyhedron $P$
with $0P = C$.

Such a polyhedron is uniquely determined by its cone at infinity $C$
and its vertices $v_i$.  Indeed, the cone over $P$ is $C(P) =
\Q_+\langle 0 \! \times \!  C, \, 1 \! \times \! v_i \rangle$, and $P
= C(P) \cap (1 \times \Q^d)$.  So it suffices to show that, at the
fixed points of the torus action on the Nash blow-up, the weights of
the torus action on $\O(1)$ are exactly the coordinates of the
vertices of $C + \Hull S$.

Our choice of an embedding $X \subset \A^n$ will be the following
canonical one.  The surjection $\Z_+^n \to C \cap \Z^d$ sending the
standard basis vectors to the rows of $M$ induces a surjection of
algebras $K[\Z_+^n] \to K[C \cap \Z^d]$.  The corresponding morphism $\Spec
K[C \cap \Z^d] \to \Spec K[\Z_+^n]$ is the desired embedding.

Let $p$ be the basepoint of $X$: the point so that for every monomial
$f \in K[C \cap \Z^d]$, $f(p)= 1$.  The homomorphisms of
algebras $$K[\Z_+^n] \to K[C \cap \Z^d] \to K[\Z^d] \to K,$$
where the last map sends every monomial to 1, correspond to the
inclusions of schemes
$$\A^n \supset X \supset T \supset \{ p \}.$$

By remark (4), we may consider the affine version of the Gauss
map for this embedding.  This is a rational map $G: X
\dasharrow \Gr(d,n)$.  We claim that $G(p)$ is the span of the
columns of $M$.  Indeed, in terms of variables $x_1, \dots,
x_n$ and $y_1, \dots, y_d$, the homomorphism $K[\Z_+^n] \to
K[\Z^d]$ is given by $x_i \mapsto \prod_j y_j^{m_{ij}}$.  The
parametric curve $y_j = 1 + t \delta_{ij}$ in $T$ therefore
maps to $x_i = (1+t)^{m_{ij}}$ in $\A^n$, so its derivative
with respect to $t$ at $0$ is $(m_{1j},\dots,m_{nj})$, the
$j$th column of $M$.

The coordinates of the Pl\"ucker embedding $\Gr(d,n) \to \P \Lambda^d
K^n$ are indexed by $d$-element subsets $I \subset \{1,\dots,n\}$.
This embedding is $T$-equivariant for the induced linear action of
$T$ on $\P \Lambda^d K^n$.  The $I$th Pl\"ucker coordinate of $G(p)$ is
the $I$th minor of $M$.  Hence $G(p)$ is contained in the linear subspace
of $\P \Lambda^d K^n$ spanned by those coordinates $I$ for which the
$I$th minor of $M$ is nonzero.  Since the $T$-action on $\P \Lambda^d
K^n$ is diagonal, the entire closure of the orbit of $G(p)$ must be
contained in this subspace.  Hence any fixed point in the closure of
this orbit must be the $I$th coordinate axis $e_I$ for some $I$ as above.
If $I = \{ i_1,\dots,i_d \}$, then the nonvanishing of the $I$th minor
is equivalent to the linear independence of $h_{i_1},\dots, h_{i_d}
\in H$, and the fiber of $\O(1)$ at this point is acted on with weight
$h_{i_1}+\cdots+ h_{i_d}$.  That is, the weights at fixed points in
this subspace are exactly the elements of $S$.

The closure of the graph of the Gauss map is clearly contained in $X
\times \overline{T p}$, so its $T$-fixed points must be of the form
$q \times e_I$, where $q$ is the unique fixed point in $X$, and $e_I$
is as above.  The weights of $\O(1)$ at these points must therefore
belong to $S$.  The same is true for the normalization, since $\O(1)$
pulls back to an ample bundle there.

Consequently, $P$ is a polyhedron with $0P = C$ and with vertices
contained in $S$.  Therefore $P \subset C + \Hull S$.

To establish equality, it suffices to show that every vertex in $C +
\Hull S$ is the weight of some fixed point in the Nash blow-up.  For
every vertex $v_I$ of $C + \Hull S$, there is a linear functional $f$
on $\Q^d$ whose restriction to $C + \Hull S$ takes on its minimum only
at $v_I$.  Hence its restriction to $\Hull S$ takes on its minimum
only at $v_I$, and its restriction to $C$ takes on its minimum only at
$0$.  The corresponding 1-parameter subgroup $\la(t): K^\times \to T$
therefore satisfies $\lim_{t \to 0} \la(t) \cdot G(p) = e_I$ and
$\lim_{t \to 0} \la(t) \cdot p = q$.  Hence $q \times e_I \in
\overline{T \cdot (p \times G(p))}$, the closure of the graph of the
Gauss map.  A point in the normalization lying over $q \times e_I$
is acted on with the same weight.  This completes the proof. \fp

\bit{Resolution trees}

We wish to consider whether a toric variety is desingularized
by a finite sequence of Nash blow-ups.  The Nash blow-up is a
local construction: that is, the Nash blow-ups of an open cover
furnish an open cover of the Nash blow-up.  Hence it suffices
to consider an affine toric variety $X(C)$.  The Nash blow-up
of $X(C)$ is $X(C + \Hull S)$; by remark (6), an open cover
of this consists of the affines $X((C + \Hull S)_v)$, where $v$
runs over the vertices of $C + \Hull S$.  By remark (8), $X(C
+ \Hull S)$ is smooth if and only if each localization $(C +
\Hull S)_v$ is equivalent to the orthant under the action of
$GL(d,\Z)$.  If not, the Nash blow-up can be repeated by
applying the theorem to each cone $(C + \Hull S)_v$.

In other words, the process of iterating Nash blow-ups of
$X(C)$ corresponds, via the toric dictionary, to the following
algorithm in convex geometry:

(1) Given the cone $C$, find the Hilbert basis $H$ of $C \cap
\Z^d$.

(2) Find $S = \{h_1 + \cdots + h_d \, | \, h_i \in H \mbox{
  linearly independent} \}$.  

(3) Find the convex hull $\Hull S$ (i.e.\ list its vertices, or
list the inequalities defining it).

(4) Find the Minkowski sum $C + \Hull S$ (i.e.\ list its
vertices and cone at infinity, or list the inequalities defining
it).

(5) Find the localization $C' = (C + \Hull S)_v$ of this
Minkowski sum at each vertex $v$.

(6) Determine whether each such $C'$ is equivalent to the orthant.  If
so, stop; if not, apply the entire algorithm to $C'$.

Because each cone may give rise to several more in step (5),
the algorithm branches.  This can be expressed in terms of a
graph as follows.  Define the {\it Nash blow-up of a cone} $C$
to be the finite set of cones of the form $(C + \Hull S)_v$,
where $S$ is as in (2), and $v$ runs over the vertices of $C +
\Hull S$.  Then define the {\it resolution tree} of $C$ or
$X(C)$ to be the unique rooted tree, with nodes labeled by
cones in $\Q^d$, whose root is labeled by $C$, and where for
every node, say labeled by $C'$:

(a) if $C'$ is equivalent to the orthant, there are no edges
beginning at $C'$ (that is, $C'$ is a leaf);

(b) otherwise, the edges beginning at $C'$ connect it to nodes
labeled by the cones $(C' + \Hull
S')_{v'}$ appearing in its Nash blow-up.

It is then clear that $X(C)$ is desingularized by a finite
number of Nash blow-ups if and only if its resolution tree is
finite.  
It is equally clear that the latter property is amenable to
computer investigation, using the algorithm above.  We will
report on this presently, but first, we explain how, in the
2-dimensional case, the situation can be completely understood.

\bit{The 2-dimensional case}

Gonzalez-Sprinberg showed \cite{cr,thesis} that toric surfaces are
desingularized by a finite sequence of (normalized) Nash
blow-ups.  This was later extended to arbitrary surfaces by
Hironaka \cite{hironaka} and Spivakovsky \cite{bourbaki,
  spivakovsky}.  In this section, we give an alternative proof of
Gonzalez-Sprinberg's original result, emphasizing the role of
Hirzebruch-Jung continued fractions.  We begin by defining them
and recalling their basic properties.

For integers $a_1, a_2, \dots$, let
$$
[a_1, \dots, a_i] =
a_1 - \cfrac{1}{a_2 - \cfrac{1}{\ddots - \cfrac{1}{a_i}}}.
$$
We assume implicitly throughout that no denominator is zero;
this is the case, for example, when $a_i > 1$ for $i > 1$.

Set $p_{-1} = 0$ and $q_0 = 0$; set $p_0 = 1$ and $q_1 = 1$.  Then
recursively let
$$
p_i = a_i p_{i-1} - p_{i-2}, \quad\quad
q_i = a_i q_{i-1} - q_{i-2}
$$
for greater values of $i$.

\bs{Proposition}
\label{proposition:hj-fraction}
For $p_i$, $q_i$ as above, $[a_1, \dots, a_i] = p_i/q_i$.
\es

\pf.  Using induction on $i$, we will prove the more general
statement where the $a_i$ are merely rational.  The case $i =
1$ is trivial. For $i>1$, assume the statement holds for
continued fractions of length $i-1$, and consider $[a_1, \dots,
a_{i-2}, a_{i-1} - 1/a_i]$.  Let $P_j, Q_j$ be the numbers
defined as above for this continued fraction.  Then $P_j = p_j$
and $Q_j = q_j$ for $j < i-1$, and
\begin{eqnarray*}
[a_1, \dots, a_i]
& = &
[a_1, \dots, a_{i-2}, a_{i-1} - 1/a_i] \\
& = &
\frac{P_{i-1}}{Q_{i-1}} \\
& = &
\frac{(a_{i-1} - 1/a_i) p_{i-2} - p_{i-3}}{(a_{i-1} - 1/a_i)q_{i-2} - q_{i-3}} \\
& = &
\frac{(a_{i-1}a_i - 1)p_{i-2} - a_i p_{i-3}}{(a_{i-1}a_i -
  1)q_{i-2} - a_i q_{i-3}} \\ 
& = &
\frac{a_i p_{i-1} - p_{i-2}}{a_i q_{i-1} - q_{i-2}} \\
& = &
\frac{p_i}{q_i}. \,\,\,\,\,\,\, \square
\end{eqnarray*}

\bs{Proposition}
\label{proposition:hj-difference}
For $i > 0$, $p_{i-1} q_i- p_i q_{i-1} = 1$.
\es

\pf.
Again use induction on $i$. The case $i = 1$ is trivial. For $i
> 1$, by the induction hypothesis, 
\begin{eqnarray*}
p_{i-1}q_i - p_i q_{i-1} & = & 
p_{i-1} (a_i q_{i-1} - q_{i-2}) - (a_i p_{i-1} - p_{i-2})
 q_{i-1} \\
& = &
-p_{i-1} q_{i-2}+ p_{i-2} q_{i-1} \\
& = & 
1. \,\,\,\,\,\,\, \square
\end{eqnarray*}

\bs{Corollary}
\label{corollary:hj-coprime}
The fraction $p_i/q_i$ is in lowest terms.  \fp
\es

\bs{Proposition}
\label{proposition:hj-partial}
For $i < j$, the denominator of $[a_{i+1}, \dots, a_j]$ as a
fraction in lowest terms is $p_i q_j - p_j q_i$.
\es

\pf.  The case $i = j-1$ is covered by the previous
proposition.  Now proceed by descending induction on $i$.  Let
$[a_{i+1}, \dots, a_j] = N_i/D_i$ in lowest terms, so that
$D_{j-1} = 1$ in particular.  Take $N_j = 1, D_j = 0$ by
convention.  Then for all $i < j$ we have
\begin{eqnarray*}
[a_{i+1}, \dots, a_j] & = & 
a_{i+1} - 1/[a_{i+2}, \dots, a_j] \\
& = &
a_{i+1} - \frac{D_{i+1}}{N_{i+1}} \\
& = &
\frac{a_{i+1} N_{i+1} - D_{i+1}}{N_{i+1}},
\end{eqnarray*} 
which is also in lowest terms.  Hence $D_i = N_{i+1}$, and
the $D_i$ satisfy the descending recursion $D_i = a_{i+2} D_{i+1} -
D_{i+2}$ with initial conditions $D_0 = 0, D_1 = 1$.  The same
holds for
\begin{eqnarray*}
p_i q_j - p_j q_i & = &
(a_{i+2} p_{i+1} - p_{i+2})q_j - p_j(a_{i+2} q_{i+1} - q_{i+2}) \\
& = &
a_{i+2}(p_{i+1}q_j - p_j q_{i+1}) - (p_{i+2}q_j - p_j q_{i+2}),
\end{eqnarray*}
which completes the proof.
\fp

\bs{Proposition}
\label{proposition:hj-algorithm}
For any rational $x$, there exists a unique finite sequence of
integers $a_1, \dots, a_k$ with $a_i > 1$ for $i > 1$ such that
$x = [a_1, \dots, a_k]$.
\es

\pf.
For any such sequence and for any $i > 1$, we have $[a_i, \dots, a_k] >
1$ by descending induction on $i$.  If $x = [a_1, \dots, a_k]$,
then $x = a_1 - 1/[a_2, \dots, a_k]$, so $a_1 =
\lceil x \rceil$ is uniquely determined.  Then $1/(a_1 - x) = [a_2, \dots,
a_k]$, and hence $a_2$ is uniquely determined too.  By
induction, all the $a_i$ are uniquely determined.

As for existence, this can be established by iterating three operations: round up,
subtract, and invert.  That is, given $x_1 = x$, let $a_1 = \lceil x_1
\rceil$, let $b_1 = a_1 - x_1$, and let $x_2 = 1/b_1$. Recursively,
given $x_i$, let $a_i = \lceil x_i \rceil$, let $b_i = a_i - x_i$,
and let $x_{i+1} = 1/b_i$.  If $x_i = n_i/d_i$ is in lowest terms,
then $x_{i+1} = d_i/(a_i d_i - n_i)$ is also in lowest terms,
so $n_{i+1} = d_i$. Since $x_i > 1$ for $i > 1$, the sequence
of $d_i$ must be nonnegative and strictly decreasing, so
eventually some $d_i = 1$ (whereupon $x_{i+1}$ is undefined and
the sequence ends).  It is then easy to verify that $x = [a_1, \dots, a_k]$.
\fp

\bs{Corollary}
\label{corollary:hj-denominators-right}
If $a_i > 1$ for $i > 1$, then for $1<i\leq j$, the denominator
of $[a_i, \dots, a_j]$ is strictly less than that of
$[a_1, \dots, a_j]$.
\es

\pf.
The sequence of denominators is the strictly decreasing
sequence $d_i$ appearing in the proof of the previous proposition.
\fp

\bs{Proposition}
\label{proposition:hj-denominators-left}
If $a_i > 1$ for $i > 1$, then for all $i < j$, the denominator
of $[a_1, \dots, a_i]$ is strictly less than that of
$[a_1, \dots, a_j]$.
\es

\pf.
The denominators are exactly the $q_i$, so this is equivalent
to showing the $q_i$ are strictly increasing, which is proved by induction on $i$:
if $q_{i-1} - q_{i-2} > 0$, then $q_i - q_{i-1} = a_i q_{i-1} -
q_{i-2} - q_{i-1} = (a_i - 1) q_{i-1} + q_{i-1} - q_{i-2} > 0$.
\fp

\bs{Corollary}
\label{corollary:hj-denominators-both}
If $a_i > 1$ for $i > 1$, then for all $1 < i < j < k$, the
denominator of $[a_i, \dots, a_j]$ is strictly less than that
of $[a_1, \dots, a_k]$.
\es

\pf.
Combine the last two results.
\fp

Now let $C$ be a proper cone in $\Q^2$. 
It can be placed in a standard form as follows.

\bs{Proposition}
\label{proposition:cone-standard}
There exists an element of $SL(2, \Z)$ taking $C$ to
$\Q_+\langle (1,0), (p, q) \rangle$ with $0 \leq p < q$ and
$p,q$ coprime; that is, a cone in the first quadrant
subtending an angle between $45^\circ$ and $90^\circ$.  
\es

\pf.
Any proper cone in $\Q^2$ has two facets or edges. Let $(a,b) \in \Z^2$
be the smallest nonzero integer point along the clockwise edge.
Then $a$ is coprime to $b$, say $ac + bd =
1$, and $\left( \begin{smallmatrix} \phantom{-}c & d \\ -b & a
\end{smallmatrix} \right) \in SL(2, \Z)$
takes $C$ to a cone whose clockwise edge is along the positive
$x$-axis and hence is contained in the first and second
quadrants. Let $(e,f)$ be the smallest nonzero integer point along the
counterclockwise edge.  Since $f > 0$, there exists
an integer $g$ such that $0 \leq e + gf < f$. Then
$\left( \begin{smallmatrix} 1 & g \\ 0 & 1 \end{smallmatrix}
\right) \in SL(2, \Z)$ takes this cone to $\Q_+\langle (1, 0),
(e+gf, f) \rangle$, which satisfies the desired properties.
\fp

In light of the last proposition, we may assume $C =
\Q_+\langle (1,0), (p,q) \rangle$ for coprime $p,q$ with $0
\leq p < q$.  As in \S2, the intersection $C \cap \Z^2$ is an
additive semigroup with a finite Hilbert basis $\H$.  In this
simple case, the Hilbert basis may be explicitly described.

\bs{Proposition}
\label{proposition:hilbert-basis}
If $p/q = [a_1, \dots, a_k]$, then $\H = \{ v_0, \dots, v_k
\}$, where $v_i = (p_i, q_i) \in \Z^2$.
\es

\pf.  Since $p_{i-1} q_i - p_i q_{i-1} = 1$, the slopes of the
rays through the $v_i$ are strictly increasing, and the
lattice points in $\Q_+\langle v_{i-1},
v_i \rangle $ are all integral linear combinations of
$v_{i-1}$ and $v_i$.  This fan of subcones covers the
entire cone, so any nonzero indecomposable element must be one
of the $v_i$.  

Conversely, since the $q_i$ are strictly increasing, if
any $v_i$ can be nontrivially expressed as an integral
combination of indecomposable elements, those elements must
belong to $\{ v_0, \dots, v_{i-1}\}$.  But this
is absurd, as those elements subtend a smaller cone that does
not contain $v_i$.  \fp

Now, as in \S\S2 and 3, let $S = \{ v_i + v_j \;|\; 0 \leq i
< j \leq k \}$.

\bs{Proposition}
\label{proposition:s-prime}
If $S' = \{v_i + v_{i+1} \;|\; 0 \leq i < k \}$, then $C +
\Hull S' = C + \Hull S$.
\es

\pf.
One inclusion is trivial.  For the other, it suffices to show that
$v_i + v_j \in C + \Hull S'$ for $0 < i+1< j \leq k$. In
fact, we will show that $v_i + v_j$ is in the even smaller set $C +
\Hull\{ v_i + v_{i+1}, v_{j-1} + v_j \}$.  This is bounded by
three lines, so it suffices to show that $v_i + v_j$ is on the
correct side of each.

First, consider the line joining $v_i + v_{i+1}$ and
$v_{j-1} + v_j$.  To simplify the notation, let $\langle (x,y), (x',
y') \rangle = x y' - x' y$, which is positive if and only if
$(x',y')$ is counterclockwise from $(x,y)$. 
By 
\re{proposition:hj-partial}, $\langle v_i, v_j \rangle >
0$ for $i < j$, and hence $\langle v_i + v_{i+1}, v_{j-1} + v_j
\rangle > 0$ too. For any two points $u_1, u_2 \in \Q^2$, the
affine linear functional $f(u) =
\langle u_1, u \rangle + \langle u, u_2 \rangle - \langle u_1,
u_2 \rangle$ vanishes on the line joining $u_1$ and
$u_2$.  Let $u_1 = v_i + v_{i+1}$ and $u_2 = v_{j-1} + v_j$; then
$f(0) = - \langle v_i + v_{i+1}, v_{j-1} + v_j \rangle < 0$.  A
brief calculation shows $f(v_i + v_j) = \langle v_i, v_j
\rangle - \langle v_i, v_{i+1} \rangle - \langle v_{i+1},
v_{j-1} \rangle - \langle v_{j-1}, v_j \rangle$.  By 
\re{proposition:hj-difference} $\langle v_i, v_{i+1} \rangle =
\langle v_{j-1}, v_j \rangle = 1$, and by 
\re{proposition:hj-partial} and
\re{corollary:hj-denominators-both} $\langle v_i, v_j \rangle
> \langle v_{i+1}, v_j \rangle > \langle v_{i+1}, v_{j-1}
\rangle$, so $\langle v_i, v_j \rangle \geq \langle v_{i+1},
v_{j-1} \rangle + 2$, and $f(v_i + v_j) \geq 0$.  Therefore
$v_i + v_j$ is on the correct side of the line.

Next, consider the line through $v_i + v_{i+1}$ with slope
$v_0$. For $v_i + v_j$ to be on the correct side of the line,
we need $(v_i + v_j) - (v_i + v_{i+1}) = v_j - v_{i+1}$ to
be counterclockwise from $v_0$, that is $\langle v_0, v_j
\rangle > \langle v_0, v_{i+1} \rangle$. Since $j \geq i+1$,
this follows from \re{proposition:hj-partial} and
\re{proposition:hj-denominators-left}.

The case of the third line is similar.
\fp

\begin{figure}

\begin{center}
\begin{picture}(245,280)

\put(123,270){$C \, \rightarrow$}
\put(165,270){$\leftarrow \, C + \Hull S$}

\thinlines

\put(0,0){\line(0,1){270}}

\multiput(20,0)(20,0){12}{\line(0,1){265}}
\multiput(0,20)(0,20){13}{\line(1,0){245}}

\Thicklines

\put(0,0){\line(1,0){250}}

\path(0,0)(160,280)

\path(18,2)(22,-2)
\path(18,-2)(22,2)

\path(18,22)(22,18)
\path(18,18)(22,22)

\path(38,62)(42,58)
\path(38,58)(42,62)

\path(58,102)(62,98)
\path(58,98)(62,102)

\path(78,142)(82,138)
\path(78,138)(82,142)

\put(40,20){\circle*{4}}
\put(60,80){\circle*{4}}
\put(100,160){\circle*{4}}
\put(140,240){\circle*{4}}

\put(60,60){\circle{4}}
\put(80,100){\circle{4}}
\put(100,140){\circle{4}}
\put(80,120){\circle{4}}
\put(100,160){\circle{4}}
\put(120,200){\circle{4}}

\put(40,20){\line(1,0){210}}
\path(40,20)(60,80)
\path(60,80)(100,160)
\path(100,160)(140,240)
\path(140,240)(163,280)

\end{picture} \bigskip \\
{\sc Figure 1.}  The case $p/q = 4/7$: $C = \Q_+\langle(1,0),(4,7)\rangle$, \\ $H = \{
\times \}$, $S = \{ \bullet, \circ \}$, $S' = \{ \bullet \}$
\end{center}

\end{figure}

So if $i$ and $j$ are not consecutive, then $v_i + v_j$ is
inessential to the shape of $C + \Hull S$.  However, if they
are consecutive, then the opposite is true, in the following sense.

\bs{Proposition}
\label{proposition:boundary}
The boundary of $C + \Hull S$ consists of the line segments joining
$v_{i-1} + v_i$ and $v_i + v_{i+1}$ for $0 < i < k$, together
with two rays starting at $v_0 + v_1$ and $v_{k-1} + v_k$ and
pointing in the directions $v_0$ and $v_k$, respectively.
\es

\pf.  Since $v_0 = (1,0)$, it suffices to show that the slopes
of these line segments are positive (or possibly $+\infty$) and
weakly decreasing.  And, finally, that the slope of $v_k$ is no
greater than the slope of the last line segment.

The line segments in question have direction $v_{i+1} -
v_{i-1}$, so it must be shown that
$$
\frac{q_{i+1} - q_{i-1}}{p_{i+1} - p_{i-1}} 
\geq \frac{q_{i+2} - q_i}{p_{i+2} - p_i}.
$$ 
Cross-multiplying and using 
\re{proposition:hj-difference} shows this to be equivalent to
$p_{i-1} q_{i+2} - p_{i+2} q_{i-1} \geq 1$, which follows from
 \re{proposition:hj-partial}.  For the last part of
the claim, it must be shown that
$$
\frac{q_k - q_{k-2}}{p_k - p_{k-2}} \geq \frac{p_k}{q_k},
$$
which follows from  \re{proposition:hj-partial},
again after cross-multiplying.  \fp

We have shown that the vertices of $C + \Hull S$ are all of the
form $v_i + v_{i+1}$.  (Although not all $v_i + v_{i+1}$ need be
vertices: see Figure 1.)  By
\re{proposition:cone-standard}, the localization
of $C$ at any such vertex can be taken by an element of
$SL(2,\Z)$ to a cone of the form $\Q_+\langle(1,0),(p',q')\rangle$
for $p,q'$ coprime and $0 \leq p' < q'$.

\bs{Proposition}
\label{proposition:cone-descend}
Unless $p = q-1$, every localized cone satisfies $q' < q$.
\es

\pf.
There are two cases: the \emph{internal case} where $0 < i <
k-1$, and the \emph{external case} where $i = 0$ or $k-1$.

In the internal case, the two edges of the
localized cone are along $(v_{i-1} + v_i) - (v_i + v_{i+1}) =
v_{i-1} - v_{i+1}$ and $(v_{i+1} + v_{i+2}) - (v_i + v_{i+1}) =
v_{i+2} - v_i$.  Hence $q' = \langle v_{i-1} - v_{i+1}, v_{i+2} -
v_i \rangle = \langle v_{i-1}, v_{i+2} \rangle - 3$ by
 \re{proposition:hj-difference}.  By \re{proposition:hj-partial}
$\langle v_{i-1}, v_{i+2} \rangle$ is the denominator
of $[a_i, a_{i+1}, a_{i+2}]$, which 
by \re{corollary:hj-denominators-both} is strictly less than $q$.

In the external case, consider first $i=0$. The two edges of
the localized cone are along $v_0$ and $(v_1 + v_2) - (v_0 + v_1) =
v_2 - v_0$, so $q' = \langle v_0, v_2 - v_0 \rangle = \langle
v_0, v_2 \rangle$, which is the denominator of $[a_1,
a_2]$.  Again by \re{corollary:hj-denominators-both}, 
this is strictly less than $q$ unless $p/q = [a_1,
a_2]$, so that $k = 2$.  If so, the condition $0 \leq p <
q$ implies $a_1 = 1$, so $p/q = (a_2-1)/a_2$ and $p=q-1$.

Likewise, when $i = k-1$, the two edges of the localized cone
are along $v_{k-2} - v_k$ and $v_k$, so $q' = \langle v_{k-2} -
v_k, v_k \rangle = \langle v_{k-2}, v_k \rangle$, which is the
denominator of $[a_{k-1}, a_k]$. Again, this is strictly less
than $q$ unless $p/q = [a_{k-1}, a_k]$, so that $k = 2$. Hence
$p = q-1$ again.  \fp

We are now in a position to prove Gonzalez-Sprinberg's result
\cite{cr,thesis}.

\bs{Theorem}
\label{theorem:toric-surface}
Any toric surface is desingularized by a finite number of Nash blow-ups.
\es

\pf.  The question is local, so it suffices to consider an
affine toric surface corresponding to a cone $\Q_+\langle
(1,0), (p,q)\rangle$, with $0 \leq p < q$ and $p,q$ coprime.
This surface is smooth if and only if $q = 1$, for only then
will the Hilbert basis consist of exactly two elements. The
previous proposition shows that $q$ is strictly decreasing
under Nash blow-ups except at external vertices for $p =
q-1$. In this case, a direct calculation shows that both
external vertices have $p'/q' = (q-2)/q$, so the denominator
will strictly decrease at the next step.  \fp

\bit{A method for enumerating simplicial cones}

In dimension $>2$, we have no general results on the resolution
of toric varieties by Nash blow-ups.  However, using a
computer, we have carried out an extensive investigation of 3-
and 4-dimensional examples.  Our primary focus is on simplicial
cones, which correspond in the toric dictionary to affine toric
orbifolds.  But, as we will see, more general cones appear in
the Nash blow-ups of simplicial cones and must be treated as
part of the recursions.

We shall begin, then, by explaining how the simplicial cones of
a given dimension $d$, or rather their equivalence classes
under the action of $GL(d,\Z)$, can be systematically
enumerated.  

Any proper cone $C \subset \Q^d$ is defined by finitely many
linear inequalities with integer coefficients, say
$\sum_{j=1}^d a_{ij}x_j \geq 0$ for $1 \leq i \leq m$.  Without
loss of generality assume that (i) no inequality is {\it
  redundant} in that it follows from the others; and (ii) for
each fixed $i$, the $a_{ij}$ are coprime.  The $m \times d$ integer
matrix $A = (a_{ij})$ is then called a {\it presentation} of
$C$.  It is unique modulo the left action of the 
group $S_m$ of permutation matrices.  To classify cones modulo
$GL(d,\Z)$, then, is equivalent to classifying integer matrices
$A$ satisfying (i) and (ii) modulo $S_m \times GL(d,\Z)$ acting
on the left and right.  This is accomplished in practice using
the following invariant.

For a cone $C$ with presentation $A$, let $\Lambda \subset \Z^d$ be
the subgroup generated by the rows of $A$.  Define the {\it index}
$I(C) \in \Z_+$ to be the index of $\Lambda$ as a subgroup of $\Z^d$.
(This is the order of the orbifold group at the
fixed point of the torus action.)  Also, if $C^* = \{ u \in \Q^d \st
\forall \, v \in C, u \cdot v \geq 0 \}$ is the dual cone, define the
{\it dual index} $I^*(C)$ to be $I(C^*)$.  Clearly $I(C)$ and $I^*(C)$
are invariant under the $GL(d,\Z)$-action.

It is, of course, nettlesome to decide whether a given matrix
satisfies the non-redundancy condition (i).  But in the
simplicial case it is easy: a presentation $A$ of a simplicial
cone is exactly a nonsingular square integer matrix satisfying
(ii).  As such, $A$ can be taken by the right action of
$GL(d,\Z)$ into {\it Hermite normal form}
\cite[4.1]{schrijver}.  This means that there exists $B \in
GL(d,\Z)$ so that $AB$ is lower triangular, with nonnegative
entries, and each row has a unique greatest entry located on
the diagonal.  Furthermore, since the entries in any given row
of $A$ are coprime, the same is true of $AB$.  These facts can
be summarized as follows.

\begin{figure}

\begin{center}
{
\label{table:dim-3}
{\footnotesize
\begin{tabular}{ l l @{\hspace{0.25in}} l l @{\hspace{0.25in}} l l @{\hspace{0.25in}} l l @{\hspace{0.25in}} l l @{\hspace{0.25in}} l l @{\hspace{0.25in}} l l}
\toprule
$I$ & $T_3(I)$ & $I$ & $T_3(I)$ & $I$ & $T_3(I)$ & $I$ & $T_3(I)$ & $I$ & $T_3(I)$ & $I$ & $T_3(I)$ & $I$ & $T_3(I)$ \\
\cmidrule(r{0.175in}){1-2}
\cmidrule(l{-0.05in}r{0.225in}){3-4}
\cmidrule(l{-0.05in}r{0.225in}){5-6}
\cmidrule(l{-0.05in}r{0.225in}){7-8}
\cmidrule(l{-0.05in}r{0.225in}){9-10}
\cmidrule(l{-0.05in}r{0.225in}){11-12}
\cmidrule(l{-0.05in}){13-14}
1  & 1   & 31 & 182 & 61 & 662  & 91  & 1679 & 121 & 2705 & 151 & 3902 & 181 & 5582  \\
2  & 2   & 32 & 227 & 62 & 693  & 92  & 1643 & 122 & 2583 & 152 & 4591 & 182 & 6595  \\
3  & 4   & 33 & 241 & 63 & 898  & 93  & 1715 & 123 & 2951 & 153 & 4872 & 183 & 6425  \\
4  & 7   & 34 & 221 & 64 & 838  & 94  & 1551 & 124 & 2919 & 154 & 4777 & 184 & 6633  \\
5  & 8   & 35 & 277 & 65 & 883  & 95  & 1825 & 125 & 3072 & 155 & 4717 & 185 & 6667  \\
6  & 11  & 36 & 311 & 66 & 915  & 96  & 2051 & 126 & 3484 & 156 & 5298 & 186 & 6729  \\
7  & 14  & 37 & 254 & 67 & 794  & 97  & 1634 & 127 & 2774 & 157 & 4214 & 187 & 6695  \\
8  & 21  & 38 & 273 & 68 & 925  & 98  & 1846 & 128 & 3211 & 158 & 4293 & 188 & 6555  \\
9  & 23  & 39 & 329 & 69 & 965  & 99  & 2110 & 129 & 3239 & 159 & 4875 & 189 & 7872  \\
10 & 25  & 40 & 381 & 70 & 1057 & 100 & 2135 & 130 & 3445 & 160 & 5555 & 190 & 7177  \\
11 & 28  & 41 & 308 & 71 & 888  & 101 & 1768 & 131 & 2948 & 161 & 5047 & 191 & 6208  \\
12 & 43  & 42 & 393 & 72 & 1206 & 102 & 2099 & 132 & 3852 & 162 & 5283 & 192 & 7942  \\
13 & 38  & 43 & 338 & 73 & 938  & 103 & 1838 & 133 & 3485 & 163 & 4538 & 193 & 6338  \\
14 & 45  & 44 & 411 & 74 & 975  & 104 & 2227 & 134 & 3105 & 164 & 5021 & 194 & 6435  \\
15 & 59  & 45 & 476 & 75 & 1254 & 105 & 2617 & 135 & 4114 & 165 & 6211 & 195 & 8569  \\
16 & 66  & 46 & 391 & 76 & 1143 & 106 & 1961 & 136 & 3709 & 166 & 4731 & 196 & 7799  \\
17 & 60  & 47 & 400 & 77 & 1219 & 107 & 1980 & 137 & 3220 & 167 & 4760 & 197 & 6600  \\
18 & 76  & 48 & 546 & 78 & 1257 & 108 & 2561 & 138 & 3763 & 168 & 6589 & 198 & 8292  \\
19 & 74  & 49 & 477 & 79 & 1094 & 109 & 2054 & 139 & 3314 & 169 & 5187 & 199 & 6734  \\
20 & 101 & 50 & 508 & 80 & 1434 & 110 & 2499 & 140 & 4454 & 170 & 5783 & 200 & 8624  \\
21 & 107 & 51 & 543 & 81 & 1350 & 111 & 2417 & 141 & 3853 & 171 & 6046 & 201 & 7727  \\
22 & 99  & 52 & 561 & 82 & 1189 & 112 & 2702 & 142 & 3479 & 172 & 5511 & 202 & 6969  \\
23 & 104 & 53 & 504 & 83 & 1204 & 113 & 2204 & 143 & 3985 & 173 & 5104 & 203 & 7933  \\
24 & 153 & 54 & 610 & 84 & 1644 & 114 & 2601 & 144 & 4668 & 174 & 5907 & 204 & 8866  \\
25 & 135 & 55 & 643 & 85 & 1473 & 115 & 2639 & 145 & 4141 & 175 & 6566 & 205 & 8153  \\
26 & 135 & 56 & 703 & 86 & 1305 & 116 & 2565 & 146 & 3675 & 176 & 6370 & 206 & 7245  \\
27 & 163 & 57 & 671 & 87 & 1507 & 117 & 2908 & 147 & 4584 & 177 & 6017 & 207 & 8762  \\
28 & 183 & 58 & 609 & 88 & 1625 & 118 & 2419 & 148 & 4113 & 178 & 5429 & 208 & 8774  \\
29 & 160 & 59 & 620 & 89 & 1380 & 119 & 2809 & 149 & 3800 & 179 & 5460 & 209 & 8311  \\
30 & 211 & 60 & 878 & 90 & 1828 & 120 & 3483 & 150 & 4894 & 180 & 7712 & 210 & 10273 \\
\bottomrule
\end{tabular} 
}
\medskip \\
{\sc Table 1.} Number of $GL(3,\Z)$-equivalence classes of
simplicial cones in 3 dimensions.
}
\end{center}

\bigskip

\begin{center}
\label{table:dim-4}
{\footnotesize
\begin{tabular}{ l l @{\hspace{0.25in}} l l @{\hspace{0.25in}} l l @{\hspace{0.25in}} l l @{\hspace{0.25in}} l l}
\toprule
$I$ & $T_4(I)$ & $I$ & $T_4(I)$ & $I$ & $T_4(I)$ & $I$ & $T_4(I)$ & $I$ & $T_4(I)$ \\
\cmidrule(r{0.225in}){1-2}
\cmidrule(l{-0.05in}r{0.225in}){3-4}
\cmidrule(l{-0.05in}r{0.225in}){5-6}
\cmidrule(l{-0.05in}r{0.225in}){7-8}
\cmidrule(l{-0.05in}){9-10}
1  & 1   & 11 & 101 & 21 & 788  & 31 & 1550 & 41 & 3399  \\
2  & 3   & 12 & 262 & 22 & 851  & 32 & 3083 & 42 & 7441  \\
3  & 7   & 13 & 154 & 23 & 682  & 33 & 2622 & 43 & 3891  \\
4  & 16  & 14 & 264 & 24 & 1778 & 34 & 2799 & 44 & 7172  \\
5  & 18  & 15 & 337 & 25 & 1037 & 35 & 2969 & 45 & 7652  \\
6  & 37  & 16 & 476 & 26 & 1338 & 36 & 5403 & 46 & 6552  \\
7  & 36  & 17 & 305 & 27 & 1530 & 37 & 2544 & 47 & 5012  \\
8  & 83  & 18 & 657 & 28 & 2123 & 38 & 3821 & 48 & 12605 \\
9  & 85  & 19 & 409 & 29 & 1288 & 39 & 4155 & 49 & 6512  \\
10 & 116 & 20 & 894 & 30 & 3006 & 40 & 6591 & 50 & 10047 \\
\bottomrule
\end{tabular}
}
\medskip \\
{\sc Table 2.} Number of $GL(4,\Z)$-equivalence classes of
simplicial cones in 4 dimensions.

\end{center}

\end{figure}

\bs{Proposition}
\label{proposition:hermite}
Every simplicial cone $C$ is equivalent to one with a presentation
$A$ which is in Hermite normal form, and each of whose rows has
coprime entries. \fp
\es

\bs{Corollary}
There are finitely many equivalence classes of simplicial cones
of dimension $d$ and index $I$.
\es

\pf.  In the simplicial case $I(C) = |\det A\,|$, so if $A$ is in
Hermite normal form, its diagonal entries multiply to $I$.
Hence there are only finitely many choices for the diagonal
entries of $A$, and so for the subdiagonal entries as well.
\fp

For a fixed value of $I$, it is now practical to enumerate the
equivalence classes of cones $C$ using
\re{proposition:hermite}.  Indeed, two matrices $A$ and $A'$
are equivalent if and only if $SAT = A'$ for some $S \in S_d$
and $T \in GL(d,\Z)$.  Detecting this is a tractable problem
for small $d$, as one can consider $A^{-1} S A'$ for all $S \in
S_d$ and see whether any of them is an integer matrix.  In this
manner, the numbers $T_d(I)$ of equivalence classes of
$d$-dimensional cones of index $I$ were determined with a
computer for small values of $I$.  These numbers are presented
in Table 1 for $d=3$ and in Table 2 for $d=4$.  A list of
explicit representatives for each of these equivalence classes,
for the first few values of $I$, is given in Table 3 for $d=3$
and in Table 4 for $d=4$.  Many cones are {\it reducible} to a
direct sum of cones of lower dimension; if so, the direct sum
in question is shown in the right-hand column of Tables 3 and
4.  By remark (3), the Nash
blow-up of a direct sum of cones is the direct sum of their
Nash blow-ups, so only irreducible cones are interesting for
our purposes.

\begin{figure}[h]
\begin{center}
{\footnotesize

\label{table:cones-dim-3}
\hspace*{-2em}
\begin{tabular}{l l l l l l l @{\hspace{0.3in}} l l l l l l l}
\toprule
Name & $I$ & $I^\ast$ & \multicolumn{3}{l}{Presentation} & Reducibility & Name & $I$ & $I^\ast$ & \multicolumn{3}{l}{Presentation} & Reducibility \\
\cmidrule(r{0.275in}){1-7}
\cmidrule(l{-0.05in}){8-14}
$C_{1,1}$ & 1 & 1  & $(e_1,$ & $e_2,$     & $e_3)$     & $A\oplus A\oplus A$& $C_{5,4}$  & 5 & 25 & $(e_1,$ & $e_2,$ & $(1,1,5))$ &                    \\
$C_{2,1}$ & 2 & 2  & $(e_1,$ & $e_2,$     & $(0,1,2))$ & $B_{2,1} \oplus A$ & $C_{5,5}$  & 5 & 25 & $(e_1,$ & $e_2,$ & $(1,2,5))$ &                    \\
$C_{2,2}$ & 2 & 4  & $(e_1,$ & $e_2,$     & $(1,1,2))$ &                    & $C_{5,6}$  & 5 & 25 & $(e_1,$ & $e_2,$ & $(2,2,5))$ &                    \\
$C_{3,1}$ & 3 & 3  & $(e_1,$ & $e_2,$     & $(0,1,3))$ & $B_{3,1} \oplus A$ & $C_{5,7}$  & 5 & 25 & $(e_1,$ & $e_2,$ & $(2,4,5))$ &                    \\
$C_{3,2}$ & 3 & 3  & $(e_1,$ & $e_2,$     & $(0,2,3))$ & $B_{3,2} \oplus A$ & $C_{5,8}$  & 5 & 25 & $(e_1,$ & $e_2,$ & $(4,4,5))$ &                    \\
$C_{3,3}$ & 3 & 9  & $(e_1,$ & $e_2,$     & $(1,1,3))$ &                    & $C_{6,1}$  & 6 & 6  & $(e_1,$ & $e_2,$ & $(0,1,6))$ & $B_{6,1} \oplus A$ \\
$C_{3,4}$ & 3 & 9  & $(e_1,$ & $e_2,$     & $(2,2,3))$ &                    & $C_{6,2}$  & 6 & 6  & $(e_1,$ & $e_2,$ & $(0,5,6))$ & $B_{6,2} \oplus A$ \\
$C_{4,1}$ & 4 & 4  & $(e_1,$ & $e_2,$     & $(0,1,4))$ & $B_{4,1} \oplus A$ & $C_{6,3}$  & 6 & 36 & $(e_1,$ & $e_2,$ & $(1,1,6))$ &                    \\
$C_{4,2}$ & 4 & 4  & $(e_1,$ & $e_2,$     & $(0,3,4))$ & $B_{4,2} \oplus A$ & $C_{6,4}$  & 6 & 18 & $(e_1,$ & $e_2,$ & $(1,2,6))$ &                    \\
$C_{4,3}$ & 4 & 16 & $(e_1,$ & $e_2,$     & $(1,1,4))$ &                    & $C_{6,5}$  & 6 & 12 & $(e_1,$ & $e_2,$ & $(1,3,6))$ &                    \\
$C_{4,4}$ & 4 & 8  & $(e_1,$ & $e_2,$     & $(1,2,4))$ &                    & $C_{6,6}$  & 6 & 6  & $(e_1,$ & $e_2,$ & $(2,3,6))$ &                    \\
$C_{4,5}$ & 4 & 8  & $(e_1,$ & $e_2,$     & $(2,3,4))$ &                    & $C_{6,7}$  & 6 & 18 & $(e_1,$ & $e_2,$ & $(2,5,6))$ &                    \\
$C_{4,6}$ & 4 & 16 & $(e_1,$ & $e_2,$     & $(3,3,4))$ &                    & $C_{6,8}$  & 6 & 6  & $(e_1,$ & $e_2,$ & $(3,4,6))$ &                    \\
$C_{4,7}$ & 4 & 2  & $(e_1,$ & $(1,2,0),$ & $(1,0,2))$ &                    & $C_{6,9}$  & 6 & 12 & $(e_1,$ & $e_2,$ & $(3,5,6))$ &                    \\
$C_{5,1}$ & 5 & 5  & $(e_1,$ & $e_2,$     & $(0,1,5))$ & $B_{5,1} \oplus A$ & $C_{6,10}$ & 6 & 18 & $(e_1,$ & $e_2,$ & $(4,5,6))$ &                    \\
$C_{5,2}$ & 5 & 5  & $(e_1,$ & $e_2,$     & $(0,2,5))$ & $B_{5,2} \oplus A$ & $C_{6,11}$ & 6 & 36 & $(e_1,$ & $e_2,$ & $(5,5,6))$ &                    \\
$C_{5,3}$ & 5 & 5  & $(e_1,$ & $e_2,$     & $(0,4,5))$ & $B_{5,3} \oplus A$ &&&&&&& \\
\bottomrule
\end{tabular}

\medskip

{\normalsize {\sc Table 3.} Classification of simplicial cones in 3 dimensions.}
}
\end{center}
\end{figure}

\bit{Results of computer investigations}

We are now in a position to describe the empirical data
obtained with a computer.  Our program, entitled {\tt resolve},
was written in the language C$++$ and relied heavily on the
Boost open-source software libraries for C$++$, especially the
linear algebra library {\tt uBLAS} of Joerg Walter and Mathias Koch
\cite{ublas}.  Our source code, as well as extensive tables of
output, are available at
$\langle$\verb+http://www.math.columbia.edu/~thaddeus/nash.html+$\rangle$.  

\begin{figure}
{\footnotesize

\begin{center}

\label{table:cones-dim-4}
\begin{tabular}{l l l l l l l l}
\toprule
Name & $I$ & $I^\ast$ & \multicolumn{4}{l}{Presentation} & Reducibility \\
\midrule
$D_{1,1}$  & 1 & 1   & $(e_1,$ & $e_2,$ & $e_3,$       & $e_4)$       & $4A$ \\
$D_{2,1}$  & 2 & 2   & $(e_1,$ & $e_2,$ & $e_3,$       & $(0,0,1,2))$ & $B_{2,1} \oplus 2A$ \\
$D_{2,2}$  & 2 & 4   & $(e_1,$ & $e_2,$ & $e_3,$       & $(0,1,1,2))$ & $C_{2,1} \oplus A$ \\
$D_{2,3}$  & 2 & 8   & $(e_1,$ & $e_2,$ & $e_3,$       & $(1,1,1,2))$ &  \\
$D_{3,1}$  & 3 & 3   & $(e_1,$ & $e_2,$ & $e_3,$       & $(0,0,1,3))$ & $B_{3,1} \oplus 2A$ \\
$D_{3,2}$  & 3 & 3   & $(e_1,$ & $e_2,$ & $e_3,$       & $(0,0,2,3))$ & $B_{3,2} \oplus 2A$ \\
$D_{3,3}$  & 3 & 9   & $(e_1,$ & $e_2,$ & $e_3,$       & $(0,1,1,3))$ & $C_{3,3} \oplus A$ \\
$D_{3,4}$  & 3 & 9   & $(e_1,$ & $e_2,$ & $e_3,$       & $(0,2,2,3))$ & $C_{3,4} \oplus A$ \\
$D_{3,5}$  & 3 & 27  & $(e_1,$ & $e_2,$ & $e_3,$       & $(1,1,1,3))$ &  \\
$D_{3,6}$  & 3 & 27  & $(e_1,$ & $e_2,$ & $e_3,$       & $(1,1,2,3))$ &  \\
$D_{3,7}$  & 3 & 27  & $(e_1,$ & $e_2,$ & $e_3,$       & $(2,2,2,3))$ &  \\
$D_{4,1}$  & 4 & 4   & $(e_1,$ & $e_2,$ & $e_3,$       & $(0,0,1,4))$ & $B_{4,1} \oplus 2A$ \\
$D_{4,2}$  & 4 & 4   & $(e_1,$ & $e_2,$ & $e_3,$       & $(0,0,3,4))$ & $B_{4,2} \oplus 2A$ \\
$D_{4,3}$  & 4 & 16  & $(e_1,$ & $e_2,$ & $e_3,$       & $(0,1,1,4))$ & $C_{4,3} \oplus A$ \\
$D_{4,4}$  & 4 & 8   & $(e_1,$ & $e_2,$ & $e_3,$       & $(0,1,2,4))$ & $C_{4,4} \oplus A$ \\
$D_{4,5}$  & 4 & 8   & $(e_1,$ & $e_2,$ & $e_3,$       & $(0,2,3,4))$ & $C_{4,5} \oplus A$ \\
$D_{4,6}$  & 4 & 16  & $(e_1,$ & $e_2,$ & $e_3,$       & $(0,3,3,4))$ & $C_{4,6} \oplus A$ \\
$D_{4,7}$  & 4 & 64  & $(e_1,$ & $e_2,$ & $e_3,$       & $(1,1,1,4))$ &  \\
$D_{4,8}$  & 4 & 32  & $(e_1,$ & $e_2,$ & $e_3,$       & $(1,1,2,4))$ &  \\
$D_{4,9}$  & 4 & 64  & $(e_1,$ & $e_2,$ & $e_3,$       & $(1,1,3,4))$ &  \\
$D_{4,10}$ & 4 & 16  & $(e_1,$ & $e_2,$ & $e_3,$       & $(1,2,2,4))$ &  \\
$D_{4,11}$ & 4 & 16  & $(e_1,$ & $e_2,$ & $e_3,$       & $(2,2,3,4))$ &  \\
$D_{4,12}$ & 4 & 32  & $(e_1,$ & $e_2,$ & $e_3,$       & $(2,3,3,4))$ &  \\
$D_{4,13}$ & 4 & 64  & $(e_1,$ & $e_2,$ & $e_3,$       & $(3,3,3,4))$ &  \\
$D_{4,14}$ & 4 & 2   & $(e_1,$ & $e_2,$ & $(0,1,2,0),$ & $(0,1,0,2))$ & $C_{4,7} \oplus A$ \\
$D_{4,15}$ & 4 & 4   & $(e_1,$ & $e_2,$ & $(0,1,2,0),$ & $(1,0,0,2))$ & $2 B_{2,1}$ \\
$D_{4,16}$ & 4 & 4   & $(e_1,$ & $e_2,$ & $(0,1,2,0),$ & $(1,1,0,2))$ &  \\
$D_{5,1}$  & 5 & 5   & $(e_1,$ & $e_2,$ & $e_3,$       & $(0,0,1,5))$ & $B_{5,1} \oplus 2A$ \\
$D_{5,2}$  & 5 & 5   & $(e_1,$ & $e_2,$ & $e_3,$       & $(0,0,2,5))$ & $B_{5,2} \oplus 2A$ \\
$D_{5,3}$  & 5 & 5   & $(e_1,$ & $e_2,$ & $e_3,$       & $(0,0,4,5))$ & $B_{5,3} \oplus 2A$ \\
$D_{5,4}$  & 5 & 25  & $(e_1,$ & $e_2,$ & $e_3,$       & $(0,1,1,5))$ & $C_{5,4} \oplus A$ \\
$D_{5,5}$  & 5 & 25  & $(e_1,$ & $e_2,$ & $e_3,$       & $(0,1,2,5))$ & $C_{5,5} \oplus A$ \\
$D_{5,6}$  & 5 & 25  & $(e_1,$ & $e_2,$ & $e_3,$       & $(0,2,2,5))$ & $C_{5,6} \oplus A$ \\
$D_{5,7}$  & 5 & 25  & $(e_1,$ & $e_2,$ & $e_3,$       & $(0,2,4,5))$ & $C_{5,7} \oplus A$ \\
$D_{5,8}$  & 5 & 25  & $(e_1,$ & $e_2,$ & $e_3,$       & $(0,4,4,5))$ & $C_{5,8} \oplus A$ \\
$D_{5,9}$  & 5 & 125 & $(e_1,$ & $e_2,$ & $e_3,$       & $(1,1,1,5))$ &  \\
$D_{5,10}$ & 5 & 125 & $(e_1,$ & $e_2,$ & $e_3,$       & $(1,1,2,5))$ &  \\
$D_{5,11}$ & 5 & 125 & $(e_1,$ & $e_2,$ & $e_3,$       & $(1,1,3,5))$ &  \\
$D_{5,12}$ & 5 & 125 & $(e_1,$ & $e_2,$ & $e_3,$       & $(1,1,4,5))$ &  \\
$D_{5,13}$ & 5 & 125 & $(e_1,$ & $e_2,$ & $e_3,$       & $(1,2,2,5))$ &  \\
$D_{5,14}$ & 5 & 125 & $(e_1,$ & $e_2,$ & $e_3,$       & $(1,2,3,5))$ &  \\
$D_{5,15}$ & 5 & 125 & $(e_1,$ & $e_2,$ & $e_3,$       & $(2,2,2,5))$ &  \\
$D_{5,16}$ & 5 & 125 & $(e_1,$ & $e_2,$ & $e_3,$       & $(2,2,4,5))$ &  \\
$D_{5,17}$ & 5 & 125 & $(e_1,$ & $e_2,$ & $e_3,$       & $(2,4,4,5))$ &  \\
$D_{5,18}$ & 5 & 125 & $(e_1,$ & $e_2,$ & $e_3,$       & $(4,4,4,5))$ &  \\
\bottomrule
\end{tabular}

\medskip

{\normalsize {\sc Table 4.} Classification of simplicial cones in 4 dimensions.}

\end{center}

}
\end{figure}

One function of the program is to enumerate the simplicial
cones of a given dimension and index, as described in the
previous section.  However, the primary function of {\tt
  resolve} is to implement the algorithm of \S3 for carrying
out the Nash blow-up and to perform it iteratively.  The C$++$
program often invokes the external programs {\tt 4ti2} \cite{4ti2},
{\tt lrs} \cite{lrs}, and {\tt qhull} \cite{qhull}, which perform isolated
parts of the computation.  Specifically, {\tt 4ti2} is used in Step 1
to find the Hilbert basis of $C \cap \Z^n$, while {\tt lrs} is used
in Steps 3 and 4 to determine the vertices of the polyhedron $C
+ \Hull S$, and the localization at each vertex; {\tt qhull} is also
used in Step 3 to simplify the determination of the convex
hull.  Because of the intensive nature of the latter
computation, Step 3 requires by far the most computing time.

We used {\tt resolve} to find {\it Nash resolutions} (that is,
finite resolution trees of Nash blow-ups) for all 1602
3-dimensional simplicial cones with $I \leq 27$ and all 201
4-dimensional simplicial cones with $I \leq 8$, following the
classification.  A few higher-dimensional cones were also
resolved, but these required considerably more time. 
To improve efficiency, {\tt resolve} ceases searching deeper in
a resolution tree whenever it reaches a simplicial cone with $I$
strictly less than the initial value,
since this has been resolved already.  However, many non-simplicial
cones are encountered in the process of resolving simplicial
cones.  So are simplicial cones with equal or greater values of
$I$.

Table 5 presents the Nash resolutions of all irreducible 3-dimensional
simplicial cones of index $I \leq 4$.  Likewise, Table 6
presents the Nash resolutions of all irreducible 4-dimensional simplicial
cones of index $I \leq 4$.  In both tables, each line
displays the rows of a presentation of a single cone.  The index
$I$ and dual index $I^*$ are shown in brackets at right.
The first line in each block of text represents the original
cone being resolved.  The singly indented lines below it show
the cones appearing in the Nash blow-up of that cone.
Subsequent to each of those, the doubly indented lines show the
cones appearing in the Nash blow-ups of those cones.

Figure 2 depicts the resolution trees of all 
3-dimensional irreducible simplicial cones of index $I \leq 6$.  Likewise,
Figure 3 depicts the resolution trees of all but one of the
4-dimensional irreducible simplicial cones of index $I \leq 5$.  (One cone
of index 5, namely $D_{5,14}$, has an enormous resolution tree
and has been omitted.)  To avoid redundancy, each tree has been
pruned of subtrees sprouting from simplicial cones that appear
elsewhere on the page.  Also, identical subtrees sprouting from
the same node have been shown only once, but with the
multiplicity appearing as a coefficient of the first cone on
the subtree.  Furthermore, a multiple branch of the form
$kC_{1,1}$ or $kD_{1,1}$ ($k$ copies of the orthant) is denoted
even more concisely by the number $k$ inside a circle.  Thus,
for example, the notation for $C_{5,4}$ is meant to convey that
a single Nash blow-up produces the 5 cones $C_{5,6}$,
$C_{3,2}$, $C_{3,2}$, $C_{1,1}$, and $C_{1,1}$.  By definition, all
leaves of a resolution tree are orthants, but this is not
immediately apparent from the diagram because of the pruning
convention just mentioned.

The cones appearing in double-outlined boxes are non-simplicial
cones, with the number of facets in parentheses.  We did not
classify these cones, so we continue their resolution trees until
they reach simplicial cones encountered before.  Evidently,
non-simplicial cones are ubiquitous even in the resolution of
simplicial cones.  (A note about the grouping of
cones by multiplicity in the figures: simplicial cones have
been grouped if and only if they are equivalent, whereas
non-simplicial cones are grouped if and only if they have
identical resolution trees.  This is a weaker condition;
in some cases, such as the $4C(4)$ in the resolution tree of
$C_{5,5}$, we know that the cones in question are not
equivalent.)

What patterns can be observed in the data?  Most obviously, all
of the thousands of cones we have studied are eventually resolved
by Nash blow-ups.  This strongly supports an affirmative
answer to Question \re{query}.

However, although the resolution seems always to exist, it also
seems to obey neither rhyme nor reason.  Almost every straightforward
conjecture one might make about patterns in the Nash resolution
seems to be false.  We have already seen, for example, that the
resolution of a simplicial cone may involve non-simplicial
cones.  One might hope that the number of facets in the cone
remains within some reasonable bound, but the resolutions of
4-dimensional simplicial cones can require cones with as many
as 10 facets, with no end in sight.  

The behavior of the indices $I$ and $I^*$ is equally
perplexing.  A 2-dimensional cone
$\Q_+\langle(1,0),(p,q)\rangle$ with $p$ coprime to $q$ has
$I(C) = I^*(C) = q$.  As we saw in
\re{proposition:cone-descend}, this is non-increasing under
Nash blow-up (indeed, decreasing except for $p$ odd and
$q=p-1$).  But $I$ and $I^*$ can increase under Nash blow-ups,
even in dimension 3 and even when the cones involved are
simplicial.  For example, $C_{6,5} = \Q_+ \langle
(1,0,0),(0,1,0),(1,3,6)\rangle$ with $I=6$ gives rise, after a
single Nash blow-up, to $\Q_+ \langle (1,3,6),(1,3,3),(2,3,6)
\rangle \cong C_{9,23}$ with $I=9$.  A glimmer of hope is
offered by $I^*$.  For, among the thousands of cones we have
examined, there appears not one example of a simplicial cone giving
rise, after a single Nash blow-up, to another simplicial cone
with greater $I^*$.  However, there are rare cases where, after
two Nash blow-ups, one obtains a simplicial cone with greater
$I^*$.  For example, $C_{9,22} =
\Q_+\langle (1,0,0),(1,3,0),(1,0,3)\rangle$ with $I^* = 3$
gives rise, after two Nash blow-ups, to $\Q_+\langle
(1,1,0),(1,0,1),(4,3,3)\rangle$ and two other cones all with
$I^*= 4$.  Moreover, there are many cases where $I^*$ increases
when one of the cones is not simplicial.  This can be seen, for
example, in the resolution tree of $C_{7,6}$, where
$\Q_+\langle (1,0,0),(0,1,0),(2,4,7),(1,1,2)\rangle$ with $I^*
= 1$ gives rise, after a single Nash blow-up, to $\Q_+\langle
(1,0,0),(0,1,0),(1,2,2)\rangle \cong C_{2,1}$ with $I^* = 2$.

The question is reminiscent of other famous iterative problems
such as the notorious Collatz conjecture \cite{lagarias}, but
in some ways it is even worse behaved.  A striking empirical
feature is the existence of simplicial cones whose Nash
resolution is vastly larger than those of other simplicial
cones with the same index. In dimension 4, for example, the
seemingly innocent $D_{5,14} = \Q_+ \langle e_1, e_2, e_3,
(1,2,3,5) \rangle$, with $I=5$, has a resolution tree with
depth 8 and 14253 cones, while no other simplicial cone with
$I=5$ needs more than depth 3 and 108 cones.  Likewise,
$D_{7,24} = \Q_+ \langle e_1, e_2, e_3, (1,2,5,7) \rangle$, with
$I=7$, has a resolution tree with depth 11 and 35299 cones,
while no other simplicial cone with $I=7$ needs more than depth
7 and 5061 cones, and only one other needs more than depth 5
and 804 cones.

In conclusion, Question \re{query} remains wide open, but we have
amassed considerable empirical evidence supporting an
affirmative answer.  In light of the 2-dimensional case, one
might hope for a proof involving some kind of
higher-dimensional analogue of continued fractions.

\newpage

{\footnotesize 
\begin{tabbing}
XX\=XXXX\=XXXX\=\=XXXXXXXXXXXXXXXXXXXXXXXXXXX\=XXXX\=XXXX\=\= \kill 
\> $C_{2,2}$: \> \> \>                           \> $C_{4,4}$: \> \> \> \\ 
\> (1,0,0),(0,1,0),(1,1,2) [2,4] \> \> \>        \> (1,0,0),(0,1,0),(1,2,4) [4,8] \> \> \> \\  
\> \>   (1,1,2),(1,0,0),(1,1,1) [1,1] \> \>      \> \>   (0,1,0),(1,2,2),(1,0,0) [2,2] \> \> \\  
\> \>   (0,1,0),(1,0,0),(1,1,1) [1,1] \> \>      \> \> \>     (0,1,0),(1,2,2),(1,1,1) [1,1] \> \\  
\> \>   (0,1,0),(1,1,2),(1,1,1) [1,1] \> \>      \> \> \>     (0,1,0),(1,0,0),(1,1,1) [1,1] \> \\  
\> \> \> \>                                      \> \>   (1,1,2),(1,2,2),(1,0,0) [2,2] \> \> \\  
\> $C_{3,3}$: \> \> \>                           \> \> \>     (1,1,2),(1,1,1),(1,2,2) [1,1] \> \\  
\> (1,0,0),(0,1,0),(1,1,3) [3,9] \> \> \>        \> \> \>     (1,0,0),(1,1,2),(1,1,1) [1,1] \> \\  
\> \>   (1,0,0),(2,2,3),(1,1,2) [1,1] \> \>      \> \>   (1,1,2),(1,2,4),(1,2,2) [2,2] \> \> \\  
\> \>   (0,1,0),(2,2,3),(1,1,2) [1,1] \> \>      \> \> \>     (1,1,2),(1,2,4),(1,2,3) [1,1] \> \\  
\> \>   (1,1,3),(1,0,0),(1,1,2) [1,1] \> \>      \> \> \>     (1,1,2),(1,2,2),(1,2,3) [1,1] \> \\  
\> \>   (0,1,0),(1,1,3),(1,1,2) [1,1] \> \>      \> \>   (1,2,4),(0,1,0),(1,2,2) [2,2] \> \> \\  
\> \>   (0,1,0),(1,0,0),(2,2,3) [3,9] \> \>      \> \> \>     (0,1,0),(1,2,4),(1,2,3) [1,1] \> \\  
\> \> \>     (1,0,0),(0,1,0),(1,1,1) [1,1] \>    \> \> \>     (0,1,0),(1,2,2),(1,2,3) [1,1] \> \\  
\> \> \>     (2,2,3),(0,1,0),(1,1,1) [1,1] \>    \> \> \> \> \\  
\> \> \>     (2,2,3),(1,0,0),(1,1,1) [1,1] \>    \> $C_{4,5}$: \> \> \> \\  
\> \> \> \>                                      \> (1,0,0),(0,1,0),(2,3,4) [4,8] \> \> \> \\  
\> $C_{3,4}$: \> \> \>                           \> \>   (0,1,0),(1,0,0),(1,1,1) [1,1] \> \> \\  
\> (1,0,0),(0,1,0),(2,2,3) [3,9] \> \> \>        \> \>   (2,3,4),(1,0,0),(1,1,1) [1,1] \> \> \\  
\> \>   (0,1,0),(1,1,1),(1,0,0) [1,1] \> \>      \> \>   (0,1,0),(1,2,2),(1,1,1) [1,1] \> \> \\  
\> \>   (2,2,3),(1,1,1),(1,0,0) [1,1] \> \>      \> \>   (2,3,4),(1,2,2),(1,1,1) [1,1] \> \> \\  
\> \>   (2,2,3),(0,1,0),(1,1,1) [1,1] \> \>      \> \> \> \> \\  
\> \> \> \>                                      \> $C_{4,6}$: \> \> \> \\  
\> $C_{4,3}$: \> \> \>                           \> (1,0,0),(0,1,0),(3,3,4) [4,16] \> \> \> \\  
\> (1,0,0),(0,1,0),(1,1,4) [4,16] \> \> \>       \> \>   (0,1,0),(1,1,1),(1,0,0) [1,1] \> \> \\  
\> \>   (1,0,0),(1,1,2),(1,1,3) [1,1] \> \>      \> \>   (3,3,4),(1,1,1),(1,0,0) [1,1] \> \> \\  
\> \>   (0,1,0),(1,1,2),(1,1,3) [1,1] \> \>      \> \>   (3,3,4),(0,1,0),(1,1,1) [1,1] \> \> \\  
\> \>   (1,1,4),(1,0,0),(1,1,3) [1,1] \> \>      \> \> \> \> \\  
\> \>   (0,1,0),(1,1,4),(1,1,3) [1,1] \> \>      \> $C_{4,7}$: \> \> \> \\  
\> \>   (0,1,0),(1,0,0),(1,1,2) [2,4] \> \>      \> (1,0,0),(1,2,0),(1,0,2) [4,2] \> \> \> \\  
\> \> \>     (1,1,2),(1,0,0),(1,1,1) [1,1] \>    \> \>   (1,0,1),(1,1,0),(1,0,0) [1,1] \> \> \\  
\> \> \>     (0,1,0),(1,0,0),(1,1,1) [1,1] \>    \> \>   (1,0,1),(1,1,0),(1,1,1) [1,1] \> \> \\  
\> \> \>     (0,1,0),(1,1,2),(1,1,1) [1,1] \>    \> \>   (1,0,2),(1,0,1),(1,1,1) [1,1] \> \> \\  
\> \> \> \>                                      \> \>   (1,2,0),(1,1,0),(1,1,1) [1,1] \> \>  
\end{tabbing} 
}

\begin{center}
{\sc Table 5.} Nash resolutions of irreducible simplicial
  cones in 3 dimensions.
\end{center}

\newpage

\bigskip\bigskip\bigskip

{\footnotesize 
\begin{tabbing}
  XXXX\=XXXX\= XXXXXXXXXXXXXXXXXXXXXXXXXX \=XXXX\=XXXX\=\=\kill
$D_{2,3}$: \> \> \>                                               $D_{4,7}$: \> \> \> \\ 
(1,0,0,0),(0,1,0,0),(0,0,1,0),(1,1,1,2) [2,8] \> \> \>            (1,0,0,0),(0,1,0,0),(0,0,1,0),(1,1,1,4) [4,64] \> \> \> \\
 \>  (0,0,1,0),(1,1,1,2),(1,0,0,0),(1,1,1,1) [1,1] \> \>          \>  (0,0,1,0),(1,0,0,0),(1,1,1,2),(1,1,1,3) [1,1] \> \> \\
 \>  (0,1,0,0),(1,1,1,2),(1,0,0,0),(1,1,1,1) [1,1] \> \>          \>  (0,1,0,0),(1,0,0,0),(1,1,1,2),(1,1,1,3) [1,1] \> \> \\
 \>  (0,1,0,0),(0,0,1,0),(1,0,0,0),(1,1,1,1) [1,1] \> \>          \>  (0,1,0,0),(0,0,1,0),(1,1,1,2),(1,1,1,3) [1,1] \> \> \\
 \>  (0,1,0,0),(0,0,1,0),(1,1,1,2),(1,1,1,1) [1,1] \> \>          \>  (0,0,1,0),(1,1,1,4),(1,0,0,0),(1,1,1,3) [1,1] \> \> \\
 \> \> \>                                                         \>  (0,1,0,0),(1,1,1,4),(1,0,0,0),(1,1,1,3) [1,1] \> \> \\
$D_{3,5}$: \> \> \>                                               \>  (0,1,0,0),(0,0,1,0),(1,1,1,4),(1,1,1,3) [1,1] \> \> \\
(1,0,0,0),(0,1,0,0),(0,0,1,0),(1,1,1,3) [3,27] \> \> \>           \>  (0,1,0,0),(0,0,1,0),(1,0,0,0),(1,1,1,2) [2,8] \> \> \\
 \>  (0,0,1,0),(1,0,0,0),(2,2,2,3),(1,1,1,2) [1,1] \> \>          \> \>    (0,0,1,0),(1,1,1,2),(1,0,0,0),(1,1,1,1) [1,1] \> \\
 \>  (0,1,0,0),(1,0,0,0),(2,2,2,3),(1,1,1,2) [1,1] \> \>          \> \>    (0,1,0,0),(1,1,1,2),(1,0,0,0),(1,1,1,1) [1,1] \> \\
 \>  (0,1,0,0),(0,0,1,0),(2,2,2,3),(1,1,1,2) [1,1] \> \>          \> \>    (0,1,0,0),(0,0,1,0),(1,0,0,0),(1,1,1,1) [1,1] \> \\
 \>  (0,0,1,0),(1,1,1,3),(1,0,0,0),(1,1,1,2) [1,1] \> \>          \> \>    (0,1,0,0),(0,0,1,0),(1,1,1,2),(1,1,1,1) [1,1] \> \\
 \>  (0,1,0,0),(1,1,1,3),(1,0,0,0),(1,1,1,2) [1,1] \> \>          \> \> \> \\
 \>  (0,1,0,0),(0,0,1,0),(1,1,1,3),(1,1,1,2) [1,1] \> \>          $D_{4,8}$: \> \> \> \\
 \>  (0,1,0,0),(0,0,1,0),(1,0,0,0),(2,2,2,3) [3,27] \> \>         (1,0,0,0),(0,1,0,0),(0,0,1,0),(1,1,2,4) [4,32] \> \> \> \\
 \> \>    (2,2,2,3),(1,0,0,0),(0,0,1,0),(1,1,1,1) [1,1] \>        \>  (1,0,0,0),(1,1,1,2),(1,1,2,2),(1,1,2,3) [1,1] \> \> \\
 \> \>    (0,1,0,0),(1,0,0,0),(0,0,1,0),(1,1,1,1) [1,1] \>        \>  (0,1,0,0),(1,1,1,2),(1,1,2,2),(1,1,2,3) [1,1] \> \> \\
 \> \>    (0,1,0,0),(2,2,2,3),(0,0,1,0),(1,1,1,1) [1,1] \>        \>  (0,0,1,0),(1,0,0,0),(1,1,2,2),(1,1,2,3) [1,1] \> \> \\
 \> \>    (0,1,0,0),(2,2,2,3),(1,0,0,0),(1,1,1,1) [1,1] \>        \>  (0,1,0,0),(0,0,1,0),(1,1,2,2),(1,1,2,3) [1,1] \> \> \\
 \> \> \>                                                         \>  (1,1,2,4),(0,0,1,0),(1,0,0,0),(1,1,2,3) [1,1] \> \> \\
$D_{3,6}$: \> \> \>                                               \>  (0,1,0,0),(1,1,2,4),(0,0,1,0),(1,1,2,3) [1,1] \> \> \\
(1,0,0,0),(0,1,0,0),(0,0,1,0),(1,1,2,3) [3,27] \> \> \>           \>  (0,1,0,0),(0,0,1,0),(1,0,0,0),(1,1,2,2) [2,4] \> \> \\
 \>  (0,0,1,0),(1,0,0,0),(2,2,3,3),(1,1,1,1) [1,1] \> \>          \> \>    (0,0,1,0),(1,1,2,2),(1,0,0,0),(1,1,1,1) [1,1] \> \\
 \>  (1,1,2,3),(1,0,0,0),(2,2,3,3),(1,1,1,1) [1,1] \> \>          \> \>    (0,1,0,0),(0,0,1,0),(1,0,0,0),(1,1,1,1) [1,1] \> \\
 \>  (0,1,0,0),(1,1,2,3),(2,2,3,3),(1,1,1,1) [1,1] \> \>          \> \>    (0,1,0,0),(0,0,1,0),(1,1,2,2),(1,1,1,1) [1,1] \> \\
 \>  (0,1,0,0),(1,1,2,3),(1,0,0,0),(1,1,1,1) [1,1] \> \>          \>  (1,1,2,4),(1,0,0,0),(1,1,1,2),(1,1,2,3) [1,1] \> \> \\
 \>  (0,1,0,0),(0,0,1,0),(2,2,3,3),(1,1,1,1) [1,1] \> \>          \>  (0,1,0,0),(1,1,2,4),(1,1,1,2),(1,1,2,3) [1,1] \> \> \\
 \>  (0,1,0,0),(0,0,1,0),(1,0,0,0),(1,1,1,1) [1,1] \> \>          \>  (0,1,0,0),(1,0,0,0),(1,1,1,2),(1,1,2,2) [2,4] \> \> \\
 \>  (0,0,1,0),(1,1,2,2),(1,0,0,0),(2,2,3,3) [1,1] \> \>          \> \>    (0,1,0,0),(1,1,1,2),(1,1,2,2),(1,1,1,1) [1,1] \> \\
 \>  (1,1,2,3),(1,1,2,2),(1,0,0,0),(2,2,3,3) [1,1] \> \>          \> \>    (1,0,0,0),(1,1,1,2),(1,1,2,2),(1,1,1,1) [1,1] \> \\
 \>  (0,1,0,0),(1,1,2,3),(1,1,2,2),(2,2,3,3) [1,1] \> \>          \> \>    (1,0,0,0),(0,1,0,0),(1,1,1,2),(1,1,1,1) [1,1] \> \\
 \>  (0,1,0,0),(0,0,1,0),(1,1,2,2),(2,2,3,3) [1,1] \> \>          \> \> \> \\
 \>  (1,1,2,3),(0,0,1,0),(1,1,2,2),(1,0,0,0) [1,1] \> \>          $D_{4,9}$: \> \> \> \\
 \>  (0,1,0,0),(1,1,2,3),(0,0,1,0),(1,1,2,2) [1,1] \> \>          (1,0,0,0),(0,1,0,0),(0,0,1,0),(1,1,3,4) [4,64] \> \> \> \\
 \> \> \>                                                         \>  (0,0,1,0),(1,0,0,0),(1,1,2,2),(1,1,3,3) [1,1] \> \> \\
$D_{3,7}$: \> \> \>                                               \>  (0,1,0,0),(0,0,1,0),(1,1,2,2),(1,1,3,3) [1,1] \> \> \\
(1,0,0,0),(0,1,0,0),(0,0,1,0),(2,2,2,3) [3,27] \> \> \>           \>  (1,1,3,4),(1,0,0,0),(1,1,2,2),(1,1,3,3) [1,1] \> \> \\
 \>  (2,2,2,3),(0,0,1,0),(1,1,1,1),(1,0,0,0) [1,1] \> \>          \>  (0,1,0,0),(1,1,3,4),(1,1,2,2),(1,1,3,3) [1,1] \> \> \\
 \>  (0,1,0,0),(0,0,1,0),(1,1,1,1),(1,0,0,0) [1,1] \> \>          \>  (1,1,3,4),(0,0,1,0),(1,0,0,0),(1,1,3,3) [1,1] \> \> \\
 \>  (0,1,0,0),(2,2,2,3),(1,1,1,1),(1,0,0,0) [1,1] \> \>          \>  (0,1,0,0),(1,1,3,4),(0,0,1,0),(1,1,3,3) [1,1] \> \> \\
 \>  (0,1,0,0),(2,2,2,3),(0,0,1,0),(1,1,1,1) [1,1] \> \>          \>  (0,0,1,0),(1,0,0,0),(1,1,1,1),(1,1,2,2) [1,1] \> \> \\
 \> \> \>                                                         \>  (0,1,0,0),(0,0,1,0),(1,1,1,1),(1,1,2,2) [1,1] \> \> \\ 
 \> \> \>                                                         \>  (1,1,3,4),(1,0,0,0),(1,1,1,1),(1,1,2,2) [1,1] \> \> \\
 \> \> \>                                                         \>  (0,1,0,0),(1,1,3,4),(1,1,1,1),(1,1,2,2) [1,1] \> \> \\
 \> \> \>                                                         \>  (0,1,0,0),(1,1,3,4),(1,0,0,0),(1,1,1,1) [1,1] \> \> \\
 \> \> \>                                                         \>  (0,1,0,0),(0,0,1,0),(1,0,0,0),(1,1,1,1) [1,1] \> \> 
\end{tabbing}
}

\begin{center}
{\sc Table 6.} Nash resolutions of irreducible simplicial
  cones in 4 dimensions.\\ 
(continued on next page)
\end{center}

\newpage

{\footnotesize 
\begin{tabbing}
XXXX\=XXXX\= XXXXXXXXXXXXXXXXXXXXXXXXXX
\=XXXX\=XXXX\=\= \kill
$D_{4,10}$: \> \> \>                                                    $D_{4,11}$: \> \> \>  \\
(1,0,0,0),(0,1,0,0),(0,0,1,0),(1,2,2,4) [4,16] \> \> \>     (1,0,0,0),(0,1,0,0),(0,0,1,0),(2,2,3,4) [4,16] \> \> \>  \\
 \>  (1,1,1,2),(0,1,0,0),(1,2,2,2),(1,0,0,0) [2,2] \> \>     \>  (1,1,1,1),(0,0,1,0),(1,1,2,2),(1,0,0,0) [1,1] \> \>  \\
 \> \>    (0,1,0,0),(1,1,1,2),(1,1,1,1),(1,2,2,2) [1,1] \>   \>  (2,2,3,4),(1,1,1,1),(1,1,2,2),(1,0,0,0) [1,1] \> \>  \\
 \> \>    (1,0,0,0),(0,1,0,0),(1,1,1,2),(1,1,1,1) [1,1] \>   \>  (0,1,0,0),(1,1,1,1),(0,0,1,0),(1,0,0,0) [1,1] \> \>  \\
 \>  (0,0,1,0),(0,1,0,0),(1,2,2,2),(1,0,0,0) [2,2] \> \>     \>  (2,2,3,4),(0,1,0,0),(1,1,1,1),(1,0,0,0) [1,1] \> \>  \\
 \> \>    (0,1,0,0),(0,0,1,0),(1,2,2,2),(1,1,1,1) [1,1] \>   \>  (0,1,0,0),(1,1,1,1),(0,0,1,0),(1,1,2,2) [1,1] \> \>  \\
 \> \>    (0,1,0,0),(0,0,1,0),(1,0,0,0),(1,1,1,1) [1,1] \>   \>  (2,2,3,4),(0,1,0,0),(1,1,1,1),(1,1,2,2) [1,1] \> \>  \\
 \>  (0,0,1,0),(1,1,1,2),(1,2,2,2),(1,0,0,0) [2,2] \> \>     \> \> \>  \\
 \> \>    (1,1,1,2),(0,0,1,0),(1,1,1,1),(1,2,2,2) [1,1] \>   $D_{4,12}$: \> \> \>  \\  
 \> \>    (1,0,0,0),(1,1,1,2),(0,0,1,0),(1,1,1,1) [1,1] \>   (1,0,0,0),(0,1,0,0),(0,0,1,0),(2,3,3,4) [4,32] \> \> \>  \\
 \>  (1,1,1,2),(1,2,2,4),(0,1,0,0),(1,2,2,2) [2,2] \> \>     \>  (0,1,0,0),(0,0,1,0),(1,2,2,2),(1,1,1,1) [1,1] \> \>  \\
 \> \>    (0,1,0,0),(1,1,1,2),(1,2,2,4),(1,2,2,3) [1,1] \>   \>  (2,3,3,4),(0,0,1,0),(1,2,2,2),(1,1,1,1) [1,1] \> \>  \\
 \> \>    (0,1,0,0),(1,1,1,2),(1,2,2,2),(1,2,2,3) [1,1] \>   \>  (1,0,0,0),(2,3,3,4),(0,0,1,0),(1,1,1,1) [1,1] \> \>  \\
 \>  (0,0,1,0),(1,2,2,4),(0,1,0,0),(1,2,2,2) [2,2] \> \>     \>  (2,3,3,4),(0,1,0,0),(1,2,2,2),(1,1,1,1) [1,1] \> \>  \\
 \> \>    (0,1,0,0),(0,0,1,0),(1,2,2,4),(1,2,2,3) [1,1] \>   \>  (1,0,0,0),(2,3,3,4),(0,1,0,0),(1,1,1,1) [1,1] \> \>  \\
 \> \>    (0,1,0,0),(0,0,1,0),(1,2,2,2),(1,2,2,3) [1,1] \>   \>  (1,0,0,0),(0,1,0,0),(0,0,1,0),(1,1,1,1) [1,1] \> \>  \\
 \>  (0,0,1,0),(1,1,1,2),(1,2,2,4),(1,2,2,2) [2,2] \> \>     \> \> \>  \\
 \> \>    (0,0,1,0),(1,1,1,2),(1,2,2,4),(1,2,2,3) [1,1] \>   $D_{4,13}$: \> \> \>  \\
 \> \>    (0,0,1,0),(1,1,1,2),(1,2,2,2),(1,2,2,3) [1,1] \>   (1,0,0,0),(0,1,0,0),(0,0,1,0),(3,3,3,4) [4,64] \> \> \>  \\
 \> \> \>                                                                      \>  (0,1,0,0),(0,0,1,0),(1,1,1,1),(1,0,0,0) [1,1] \> \>  \\
 \> \> \>                                                                      \>  (3,3,3,4),(0,0,1,0),(1,1,1,1),(1,0,0,0) [1,1] \> \>  \\
 \> \> \>                                                                      \>  (3,3,3,4),(0,1,0,0),(1,1,1,1),(1,0,0,0) [1,1] \> \>  \\
 \> \> \>                                                                      \>  (3,3,3,4),(0,1,0,0),(0,0,1,0),(1,1,1,1) [1,1] \> \>                      
\end{tabbing}
}

\begin{center}
{\sc Table 6.} Nash resolutions of irreducible simplicial
  cones in 4 dimensions.\\
(continued from previous page)
\end{center}

\begin{center}
\label{diagram:trees-dim-3}
\begin{tikz-trees-diagram}


  \begin{scope}[xshift=0in, yshift=0in]
    \node[scone] {$C_{2,2}$}
      child{ node[leaf] {3} };
  \end{scope}


  \begin{scope}[xshift=1in, yshift=0in]
    \node[scone] {$C_{3,3}$}
      child{ node[scone] {$C_{3,4}$} }
      child{ node[leaf] {4} };
  \end{scope}

  \begin{scope}[xshift=2in, yshift=0in]
    \node[scone] {$C_{3,4}$}
      child{ node[leaf] {3} };
  \end{scope}


  \begin{scope}[xshift=3in, yshift=0in]
    \node[scone] {$C_{4,3}$}
      child { node[scone] {$C_{2,2}$} }
      child { node[leaf] {4} };
  \end{scope}

  \begin{scope}[xshift=4in, yshift=0in]
    \node[scone] {$C_{4,4}$}
      child { node[scone] {$4\,C_{2,1}$} };
  \end{scope}

  \begin{scope}[xshift=5in, yshift=0in]
    \node[scone] {$C_{4,5}$}
      child { node[leaf] {4} };
  \end{scope}

  \begin{scope}[xshift=0in, yshift=-1in]
    \node[scone] {$C_{4,6}$}
      child { node[leaf] {3} };
  \end{scope}

  \begin{scope}[xshift=1in, yshift=-1in]
    \node[scone] {$C_{4,7}$}
      child { node[leaf] {4} };
  \end{scope}


  \begin{scope}[xshift=0.5in, yshift=-2in]
    \node[scone] {$C_{5,4}$}
      child { node[scone] {$C_{5,6}$} }
      child { node[scone] {$2\,C_{3,2}$} }
      child { node[leaf] {2} };
  \end{scope}

  \begin{scope}[xshift=3.7in, yshift=-1in]
    \node[scone] {$C_{5,5}$}
      child { node at (-1,0) [nscone] {$C(5)$}
        child { node at (-0.4,0) [nscone] {$C(6)$}
          child { node[leaf] {6} }
        }
        child { node[nscone] {$4\,C(4)$}
          child { node[leaf] {4} }
        }
        child { node[leaf] {4} }
      }
      child { node at (1,0) [nscone] {$C(4)$}
        child { node[nscone] {$2\,C(4)$}
          child { node[leaf] {4} }
        }
        child { node[leaf] {5} }
      }
      child { node[scone] at (1.3,0) {$2\,C_{2,1}$} }
      child { node[leaf] at (1.3,0) {4} };
  \end{scope}

  \begin{scope}[xshift=0in, yshift=-3in]
    \node[scone] {$C_{5,6}$}
      child { node[scone] {$C_{3,4}$} }
      child { node[leaf] {4} };
  \end{scope}

  \begin{scope}[xshift=1in, yshift=-3in]
    \node[scone] {$C_{5,7}$}
      child { node[leaf] {7} };
  \end{scope}

  \begin{scope}[xshift=2in, yshift=-3in]
    \node[scone] {$C_{5,8}$}
      child { node[leaf] {3} };
  \end{scope}


  \begin{scope}[xshift=3.5in, yshift=-3in]
    \node[scone] {$C_{6,3}$}
      child { node[scone] {$C_{3,3}$} }
      child { node[scone] {$2\,C_{2,1}$} }
      child { node[leaf] {2} };
  \end{scope}

  \begin{scope}[xshift=1in, yshift=-4in]
    \node[scone] {$C_{6,4}$}
      child { node[nscone] at (-0.3,0) {$4\,C(4)$}
        child { node[leaf] {4} }
      }
      child { node[scone] {$C_{3,3}$} }
      child { node[scone] {$2\,C_{2,1}$} }
      child { node[leaf] {1} };
  \end{scope}

  \begin{scope}[xshift=2.8in, yshift=-4in]
    \node[scone] {$C_{6,5}$}
      child { node[scone] at (-0.1,0) {$2\,C_{9,23}$}
        child { node[nscone] {$2\,C(4)$}
          child { node[nscone] {$C(4)$}
            child { node[leaf] {4} }
          }
          child { node[leaf] {3} }
        }
        child { node[leaf] {2} }
      }
      child { node[scone] at (0.1,0)  {$2\,C_{3,1}$} };
  \end{scope}

  \begin{scope}[xshift=4.7in, yshift=-4in]
    \node[scone] {$C_{6,6}$}
      child { node[nscone] at (-0.5,0) {$C(4)$}
        child { node[leaf] {4} }
      }
      child { node[scone] {$2\,C_{3,2}$} }
      child { node[leaf] {3} };
  \end{scope}

  \begin{scope}[xshift=4in, yshift=-5.5in]
    \node[scone] {$C_{6,7}$}
      child { node[scone] {$2\,C_{3,2}$} }
      child { node[leaf] {4} };
  \end{scope}

  \begin{scope}[xshift=1in, yshift=-5.5in]
    \node[scone] {$C_{6,8}$}
      child { node[nscone] at (-0.3,0) {$C(4)$}
        child { node[nscone] {$C(4)$}
          child { node[leaf] {4} }
        }
        child { node[leaf] {3} }
      }
      child { node[scone] {$C_{3,2}$} }
      child { node[leaf] {2} };
  \end{scope}

  \begin{scope}[xshift=5in, yshift=-5.5in]
    \node[scone] {$C_{6,9}$}
      child { node[leaf] {4} };
  \end{scope}

  \begin{scope}[xshift=4in, yshift=-6.5in]
    \node[scone] {$C_{6,10}$}
      child { node[leaf] {4} };
  \end{scope}

  \begin{scope}[xshift=5in, yshift=-6.5in]
    \node[scone] {$C_{6,11}$}
      child { node[leaf] {3} };
  \end{scope}

\end{tikz-trees-diagram}
\end{center}

\begin{center}
{\sc Figure 2.}  Resolution trees of irreducible
simplicial cones in 3 dimensions.
\end{center}

\begin{center}
\label{diagram:trees-dim-4}
\begin{tikz-trees-diagram}


  \begin{scope}[xshift=0in, yshift=0in]
    \node[scone] {$D_{2,3}$}
      child{ node[leaf] {4} };
  \end{scope}


  \begin{scope}[xshift=1in, yshift=0in]
    \node[scone] {$D_{3,5}$}
      child{ node[scone] {$D_{3,7}$} }
      child{ node[leaf] {6} };
  \end{scope}

  \begin{scope}[xshift=2in, yshift=0in]
    \node[scone] {$D_{3,6}$}
      child{ node[leaf] {12} };
  \end{scope}

  \begin{scope}[xshift=3in, yshift=0in]
    \node[scone] {$D_{3,7}$}
      child{ node[leaf] {4} };
  \end{scope}


  \begin{scope}[xshift=4in, yshift=0in]
    \node[scone] {$D_{4,7}$}
      child{ node[scone] {$D_{2,3}$} }
      child{ node[leaf] {6} };
  \end{scope}

  \begin{scope}[xshift=5in, yshift=0in]
    \node[scone] {$D_{4,8}$}
      child{ node[scone] {$2\,D_{2,2}$} }
      child{ node[leaf] {8} };
  \end{scope}

  \begin{scope}[xshift=0in, yshift=-1in]
    \node[scone] {$D_{4,9}$}
      child{ node[leaf] {12} };
  \end{scope}

  \begin{scope}[xshift=1in, yshift=-1in]
    \node[scone] {$D_{4,10}$}
      child{ node[scone] {$6\,D_{2,1}$} };
  \end{scope}

  \begin{scope}[xshift=2in, yshift=-1in]
    \node[scone] {$D_{4,11}$}
      child{ node[leaf] {6} };
  \end{scope}

  \begin{scope}[xshift=3in, yshift=-1in]
    \node[scone] {$D_{4,12}$}
      child{ node[leaf] {6} };
  \end{scope}

  \begin{scope}[xshift=4in, yshift=-1in]
    \node[scone] {$D_{4,13}$}
      child{ node[leaf] {4} };
  \end{scope}

  \begin{scope}[xshift=5in, yshift=-1in]
    \node[scone] {$D_{4,16}$}
      child{ node[leaf] {7} };
  \end{scope}


  \begin{scope}[xshift=0.5in, yshift=-2in]
    \node[scone] {$D_{5,9}$}
      child{ node[scone] {$D_{5,15}$} }
      child{ node[scone] {$3\,D_{3,2}$} }
      child{ node[leaf] {3} };
  \end{scope}

  \begin{scope}[xshift=4.5in, yshift=-2in]
    \node[scone] {$D_{5,10}$}
      child{ node at (-1,0) [nscone] {$2\,D(6)$}
        child{ node at (-0.4,0) [nscone] {$D(7)$}
          child{ node[leaf] {6} }
        }
        child{ node[nscone] {$4\,D(5)$}
          child{ node[leaf] {4} }
        }
        child{ node[leaf] {4} }
      }
      child{ node at (1,0) [nscone] {$D(5)$}
        child{ node[nscone] {$4\,D(5)$}
          child{ node[leaf] {4} }
        }
        child{ node[leaf] {8} }
      }
      child{ node at (1.2,0) [leaf] {14} };
  \end{scope}

  \begin{scope}[xshift=1.5in, yshift=-4in]
    \node[scone] {$D_{5,11}$}
      child{ node at (-1,0) [nscone] {$2\,D(6)$}
        child{ node at (-0.4,0) [nscone] {$D(7)$}
          child{ node[leaf] {6} }
        }
        child{ node[nscone] {$4\,D(5)$}
          child{ node[leaf] {4} }
        }
        child{ node[leaf] {4} }
      }
      child{ node at (0,0) [nscone] {$2\,D(5)$}
        child{ node[leaf] {5} }
      }
      child{ node at (0.4,0) [scone] {$2\,D_{2,2}$} }
      child{ node at (0.4,0) [leaf] {10} };
  \end{scope}

  \begin{scope}[xshift=2in, yshift=-2in]
    \node[scone] {$D_{5,12}$}
      child{ node[scone] {$8\,D_{3,2}$} }
      child{ node[leaf] {4} };
  \end{scope}

  \begin{scope}[xshift=4.5in, yshift=-4in]
    \node[scone] {$D_{5,13}$}
      child{ node at (-0.8,0) [nscone] {$2\,D(5)$}
        child{ node[nscone] {$2\,D(5)$}
          child{ node[leaf] {4} }
        }
        child{ node[leaf] {1} }
      }
      child{ node at (0.5,0) [nscone] {$D(5)$}
        child{ node[nscone] {$4\,D(5)$}
          child{ node[leaf] {4} }
        }
        child{ node[leaf] {6} }
      }
      child{ node at (0.8,0) [scone] {$3\,D_{2,1}$} }
      child{ node at (0.8,0) [leaf] {9} };
  \end{scope}

  \begin{scope}[xshift=2.3in, yshift=-5in]
    \node[scone] {$D_{5,14}$}
      child{ node {$\cdots$} };
  \end{scope}

  \begin{scope}[xshift=0in, yshift=-3in]
    \node[scone] {$D_{5,15}$}
      child{ node[scone] {$D_{3,7}$} }
      child{ node[leaf] {6} };
  \end{scope}

  \begin{scope}[xshift=1in, yshift=-3in]
    \node[scone] {$D_{5,16}$}
      child{ node[leaf] {12} };
  \end{scope}

  \begin{scope}[xshift=2in, yshift=-3in]
    \node[scone] {$D_{5,17}$}
      child{ node[leaf] {10} };
  \end{scope}

  \begin{scope}[xshift=5.5in, yshift=-5in]
    \node[scone] {$D_{5,18}$}
      child{ node[leaf] {4} };
  \end{scope}

\end{tikz-trees-diagram}
\end{center}

\begin{center}
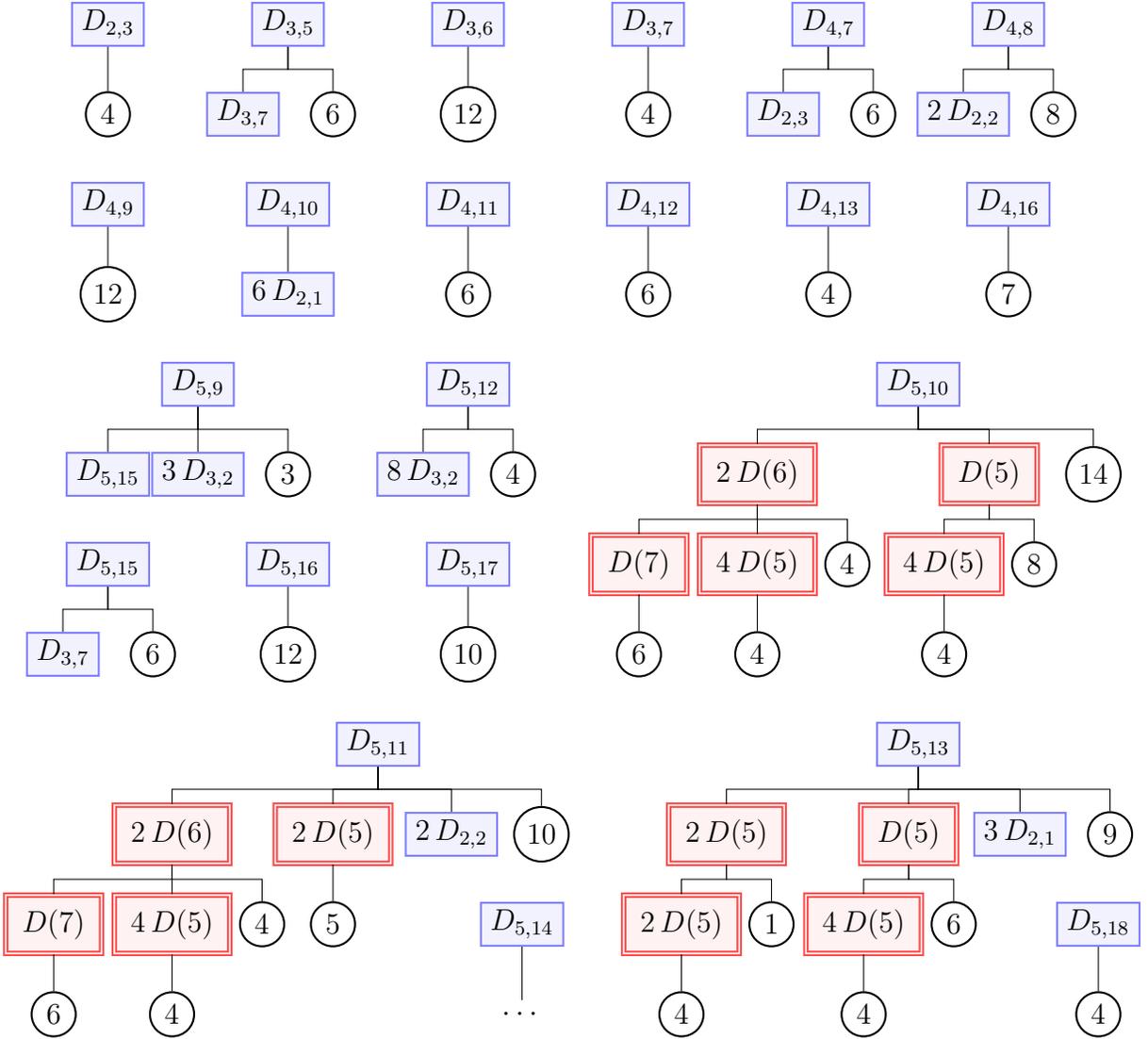
{\sc Figure 3.}  Resolution trees of irreducible
simplicial cones in 4 dimensions.
\end{center}

\end{document}